\documentclass[preprint,sort&compress]{elsarticle}
\usepackage[margin=1in]{geometry}
\usepackage{lineno}
\usepackage{amssymb,amsxtra,amsmath}
\usepackage{latexsym}
\usepackage{ifthen}
\usepackage{graphicx}
\usepackage{calc,fancyhdr}
\usepackage{indentfirst}
\usepackage{parskip,algorithm,algorithmicx}
\usepackage{algpseudocode}
\usepackage{color}
\usepackage{bigints}
\usepackage{subfig}
\usepackage{booktabs}
\usepackage{Article}
\usepackage{color}
\newcommand{\revone}[1]{\color{black}#1\normalcolor}
\newcommand{\revtwo}[1]{\color{black}#1\normalcolor}

\journal{JCP}

\begin{document}

\begin{frontmatter}

\title{Hyperviscosity-Based Stabilization for Radial Basis Function-Finite Difference (RBF-FD) Discretizations of Advection-Diffusion Equations}

\author[addr1]{Varun Shankar\corref{corresp}}
\address[addr1]{Department of Mathematics and School of Computing, University of Utah, UT, USA}
\ead{vshankar@math.utah.edu}
\cortext[corresp]{Corresponding Author}

\author[addr2]{Aaron L. Fogelson}
\address[addr2]{Departments of Mathematics and Bioengineering, University of Utah, UT, USA}
\ead{fogelson@math.utah.edu}

\begin{abstract}
We present a novel hyperviscosity formulation for stabilizing RBF-FD discretizations of the advection-diffusion equation. The amount of hyperviscosity is determined quasi-analytically for commonly-used explicit, implicit, and implicit-explicit (IMEX) time integrators by using a simple 1D semi-discrete Von Neumann analysis. The analysis is applied to an analytical model of spurious growth in RBF-FD solutions that uses auxiliary differential operators mimicking the undesirable properties of RBF-FD differentiation matrices. The resulting hyperviscosity formulation is a generalization of existing ones in the literature, but is free of any tuning parameters and can be computed efficiently. To further improve robustness, we introduce a simple new scaling law for polynomial-augmented RBF-FD that relates the degree of polyharmonic spline (PHS) RBFs to the degree of the appended polynomial. When used in a novel ghost node formulation in conjunction with the recently-developed overlapped RBF-FD method, the resulting method is robust and free of stagnation errors. We validate the high-order convergence rates of our method on 2D and 3D test cases over a wide range of Peclet numbers (1-1000). We then use our method to solve a 3D coupled problem motivated by models of platelet aggregation and coagulation, again demonstrating high-order convergence rates.
\end{abstract}

\end{frontmatter}

\section{Introduction}
\label{sec:intro}

Collocation methods based on radial basis functions (RBFs) are increasingly popular, due to their high-order convergence rates and their ability to naturally handle scattered node layouts on arbitrary domains. RBF interpolants can be used to generate both pseudospectral (RBF-PS) and finite-difference (RBF-FD) methods~\cite{Bayona2010,Davydov2011,Wright200699,FlyerNS,FlyerPHS,BarnettPHS}. RBF-based methods are also easily applied to the solution of PDEs on node sets that are not unisolvent for polynomials, such as ones lying on the sphere $\mathbb{S}^2$~\cite{FlyerWright:2007,FlyerWright:2009,FoL11,FlyerLehto2012} and other general surfaces~\cite{Piret2012, Piret2016,FuselierWright2013,SWFKJSC2014,LSWSISC2017}.

Despite their strengths, RBF methods have historically had a major drawback: poor robustness/stability. There are two major aspects to the stability of RBF-FD methods for PDEs: \emph{local} stability, and \emph{global} stability. The former refers to the stability of the process of generating RBF-FD weights (interpolation), and the latter to the eigenspectrum of the sparse matrix containing the RBF-FD weights. With regards to local stability, traditional RBF interpolation (using a shape parameter) suffered from ill-conditioning, causing researchers to posit an ``uncertainty principle'' (a tradeoff between stability and accuracy); for a review, see~\cite{Fasshauer:2007}. Fortunately, this ``local stability'' problem has been solved for RBFs with shape parameters~\cite{FoWr,FornbergPiret:2007,FaMC12,FLF,FoLePo13,FoWr2016}, and polyharmonic spline (PHS) RBFs~\cite{FlyerNS,BarnettPHS,FlyerPHS,FlyerElliptic}. In the case of RBF-FD methods based on PHS RBFs, the key to ameliorating stability issues is the inclusion of polynomials in the approximation space and enforcing polynomial reproduction on the RBF-FD weights. \revone{Typically, however, the stencil size is required to be twice the number of polynomial basis functions~\cite{FlyerElliptic,FlyerNS}. Thus, this technique results in larger stencils than if one were to use RBFs or polynomials alone}. Fortunately, the cost associated with the increase in stencil sizes for this new \emph{augmented RBF-FD} method can be dramatically reduced by using a variant called the \emph{overlapped RBF-FD} method~\cite{ShankarJCP2017}, developed by the first author and used in this article. \revone{Interestingly, there is also evidence that the inclusion of polynomials improves the global stability of RBF-FD methods, at least in the context of elliptic PDEs~\cite{FlyerElliptic}. Unfortunately, in our experiments with parabolic, hyperbolic, and mixed-character PDEs, we nevertheless encountered spurious eigenvalues in RBF-FD differentiation matrices corresponding to those differential operators}. Specifically, RBF-FD differentiation matrices may contain rogue eigenvalues with positive real parts, resulting in spurious growth in the numerical solution to time-dependent PDEs. \revone{In addition, when used with implicit time-stepping, we observed that the time-stepping matrices had spectral radii greater than unity.}

While the theory of global stability of RBF-FD methods is virtually non-existent, Fornberg and Lehto~\cite{FoL11} devised a practical approach to rectify the spectra of RBF-FD differentiation matrices corresponding to hyperbolic operators on the sphere by the addition of artificial hyperviscosity. They demonstrated that the technique successfully shifted any rogue eigenvalues to the left half of the complex plane. The hyperviscosity operator takes the form $\gamma \Delta^k$, where $\gamma \in \mathbb{R}$ and $k \in \mathbb{N}$ are small numbers that must be selected. Then, in~\cite{FlyerLehto2012}, Flyer et al. empirically computed $\gamma$ and $k$ for the shallow water equations on the sphere with two specific prescriptions: that the hyperviscosity operator vanish under refinement so the correct solution is recovered, and that the order of the operator $k$ must be increased with the stencil size $n$; this latter condition is seen in spectral methods~\cite{MaSSV2}. Recently, in~\cite{FlyerNS}, Flyer et al. use $\gamma = (-1)^{k+1} 2^{-6} h^{2k}$; however, unlike in~\cite{FlyerLehto2012}, the authors of~\cite{FlyerNS} do not increase $k$ with $n$. In his PhD thesis~\cite{BarnettPHS}, Barnett sets $\gamma = (-1)^{k+1} 2^{-6} h^{2k+1}$, also resulting in stable simulations in Euclidean domains and on the sphere. In both these cases, $2^{-6}$ is a user-defined parameter whose origin is not explained. Further, and somewhat surprisingly, these expressions are independent of the magnitude of the velocity field in which the quantity of interest is advected. It is worth noting that hyperviscosity-based stabilization also has a long history of use in spectral methods (for example,~\cite{Tadmor89,MaSSV1,MaSSV2}), but in that literature, the amount of hyperviscosity is typically determined analytically.

RBF-FD methods with hyperviscosity appear to be stable despite these inconsistencies in the value of $\gamma$, possibly due to the fact that~\cite{FoL11,FlyerLehto2012,FlyerNS,BarnettPHS} all use hand-tuned parameters in the formula for $\gamma$. Unfortunately, when we used the above recipes for hyperviscosity in the context of implicit-explicit (IMEX) linear multistep methods for solving advection-diffusion equations
\begin{align}
\frac{\partial c}{\partial t} + \vu \cdot \nabla c = \nu \Delta c,
\label{eq:ad}
\end{align}
we observed mild to severe instabilities over a range of Peclet numbers despite the presence of natural diffusion in the problem, especially on scattered and coarse node sets. 

The primary contribution of this work is to carefully derive a \emph{quasi-analytic} expression for $\gamma$, which in turn allows us to develop high-order numerical discretizations of advection-diffusion equations on scattered node sets and irregular domains. Our new expression for $\gamma$ gives the correct power of $h$ to use, an explanation for the $2^{-6}$ value, and the correct dependence on the velocity $\vu$ and the diffusion coefficient $\nu$. The resulting expression is free of any hand-tuned parameters, generalizes existing formulas in the RBF-FD literature, and connects the RBF-FD hyperviscosity literature to the spectral methods literature by which it was inspired. To derive this expression, we develop a novel semi-discrete Von Neumann analysis technique based on \emph{auxiliary} differential operators that explicitly model the spurious growth modes introduced by RBF-FD discretizations of elliptic and hyperbolic differential operators. While we focus on advection-diffusion equations, this analysis carries over naturally to any linear PDE (or linearized nonlinear PDEs). For completeness, we also present a scaling law that relates $k$ and the stencil size $n$. A secondary  contribution of this work is a new scaling law connecting the power $m$ of PHS RBFs and the degree $\ell$ of the appended polynomial; we find that this law significantly improves stability over fixing $m$ as $\ell$ is increased, or scaling $m$ as $2\ell+1$. A third contribution of this work is a new ghost node method for IMEX time-stepping of advection-diffusion equations, with both diffusion and hyperviscosity being handled \emph{implicitly} in time. \revone{While ghost nodes are not always required with RBF-FD methods~\cite{FlyerElliptic}, we found in our experiments that they help ameliorate global stability issues when used in conjunction with hyperviscosity.}

The remainder of the paper is organized as follows. In Section 2, we review the overlapped RBF-FD method (with error estimates). In Section 3, we present our new scaling law for relating the PHS degree $m$ to the polynomial degree $\ell$, and compare this to other possible laws. Then, in Section 4, we derive quasi-analytic expressions for the hyperviscosity parameter $\gamma$ for different time integrators (explicit, implicit, and IMEX), and present efficient techniques for computing spurious growth modes; we also present an expression for $k$, and discuss how to approximate the hyperviscosity operator efficiently with overlapped RBF-FD. Section 5 contains our new ghost node formulation for IMEX time-stepping of advection-diffusion equations; we also briefly discuss how to deal with filling ghost values at the initial time-step. We present 2D and 3D convergence studies in Section 6, showing results over three orders of magnitude of Peclet numbers. We apply our new techniques to solving a coupled problem inspired by mathematical models of platelet aggregation and coagulation in Section 7, and show high orders of convergence on problem in the spherical shell. We conclude with a summary and comments on future work in Section 8.
\section{Overlapped RBF-FD}
\label{sec:overlap}

\subsection{Mathematical description}
\label{sec:orbf_math}

We first present a mathematical description of the overlapped RBF-FD method, recently developed by the first author~\cite{ShankarJCP2017}. Let $X = \{\vx_k\}_{k=1}^N$ be a global set of nodes on the domain. Define the stencil $P_k$ to be the set of nodes containing node $\vx_{\calI^k_1}$ and its $n-1$ nearest neighbors $\{\vx_{\calI^k_2},\hdots,\vx_{\calI^k_n}\}$; here, $\{\calI^k_1,\hdots,\calI^k_n\}$ are indices that map into the global node set $X$. We defer discussion of the number of stencils to the end of this section. For the remainder of this discussion, we will focus without loss of generality on the stencil $P_1$. First, define the \emph{stencil width} $\rho_1$ as
\begin{align}
\rho_1 = \max\limits_{j} \|\vx_{\calI^1_1} - \vx_{\calI^1_j}\|, j=1,\hdots,n.
\label{eq:stencil_rad}
\end{align}
Choosing $\delta \in (0,1]$, we now define the \emph{stencil retention distance} $r_1$ to be
\begin{align}
r_1 = (1-\delta)\rho_1.
\label{eq:ball_rad}
\end{align}
Here, $\delta$ is called the \emph{overlap parameter}. Let $R_1$ be the set of global indices of the $p_1$ nodes in the stencil $P_1$ that are within the distance $r_1$ from $\vx_{\calI^1_1}$. We write $R_1$ as
\begin{align}
R_1 = \{\calR^1_1, \calR^1_2, \hdots, \calR^1_{p_1}\}.
\label{eq:ball_inds}
\end{align}
In general, $R_1$ is some permutation of a subset of the indices of the nodes in $P_1$. Let $\mathbb{B}_1$ contain the nodes whose indices are in $R_1$. Thus, 
\begin{align}
\mathbb{B}_1 = \{\vx_{\calR^1_1}, \hdots, \vx_{\calR^1_{p_1}} \}.
\label{eq:ball_pts}
\end{align}
For convenience, we will refer to $\mathbb{B}_1$ as a \emph{ball}, though it is merely a set of discrete points. The overlapped RBF-FD method involves computing RBF-FD weights for all the nodes in the ball $\mathbb{B}_1$, and repeating this process for each stencil $P_k$. The weights for all the nodes in $\mathbb{B}_1$ with indices in $R_1$ are computed using the following augmented local RBF interpolant on $P_1$:
\begin{align}
s_1(\vx,\vy) = \sum\limits_{j=1}^n w^1_j(\vy) \|\vx - \vx_{\calI^1_j}\|^m + \sum\limits_{i=1}^{M} \lambda^1_i(\vy) \psi^1_i(\vx),
\label{eq:rbf_interp}
\end{align}
where $\|\vx - \vx_{\calI^1_j}\|^m$ is the polyharmonic spline (PHS) RBF of degree $m$ ($m$ is odd), and $\psi^1_i(\vx)$ are the $M$ monomials corresponding to the total degree polynomial of degree $\ell$ in $d$ dimensions. Here, the $n$ overlapped RBF-FD weights associated with the point $\vy$ are $w^1_j(\vy), j=1,\hdots,n$. \revone{Thus, for the stencil $P_1$, we have two sets of points that are relevant to the interpolant: the set $P_1$ itself (with $\vx \in P_1$), and the set $\mathbb{B}_1$ (with $\vy \in \mathbb{B}_1$). In order to compute the weights for the linear operator $\calL$ uniquely at all nodes in $\mathbb{B}_1$ with indices in the set $R_1$, we impose the following two (sets of) conditions:
\begin{align}
\lf.s_1 \rt|_{\vx \in P_1, \vy \in \mathbb{B}_1} &= \lf.\calL \|\vx - \vx_{\calI^1_j}\|^m\rt|_{\vx \in \mathbb{B}_1}, j=1,\hdots,n, \label{eq:interp_constraint}\\
\sum_{j=1}^n \lf.w_j^1(\vy) \psi_i^1(\vx) \rt|_{\vx \in P_1, \vy \in \mathbb{B}_1} &= \lf.\calL \psi^1_i(\vx)\rt|_{\vx \in \mathbb{B}_1}, i=1,\hdots,M. \label{eq:poly_constraint}
\end{align}
}
The first set of conditions enforces that $s_1(\vx,\vy)$ interpolate the derivatives of the PHS RBF at all nodes in $\mathbb{B}_1$. The second set of conditions enforces polynomial reproduction on the overlapped RBF-FD weights. If $\ell$ is the degree of the appended polynomial $\psi^1(\vx)$, $M = {\ell + d \choose d}$; for stability, we also require that $M \leq \lfloor \frac{n}{2} \rfloor$~\cite{ShankarJCP2017,FlyerNS,FlyerPHS}. The interpolant \eqref{eq:rbf_interp} and the two conditions \eqref{eq:interp_constraint}--\eqref{eq:poly_constraint} can be collected into the following block linear system:
\begin{align}
\begin{bmatrix}
A_1 & \Psi_1 \\
\Psi_1^T & O
\end{bmatrix}
\begin{bmatrix}
W_1 \\
\Lambda_1
\end{bmatrix}
=
\begin{bmatrix}
B_{A_1} \\
B_{\Psi_1}
\end{bmatrix},
\label{eq:rbf_linsys}
\end{align}
where
\begin{align}
(A_1)_{ij} &= \|\vx_{\calI^1_i} - \vx_{\calI^1_j} \|^m, i,j=1,\hdots,n, \\
(\Psi_1)_{ij} &= \psi^1_j(\vx_{\calI^1_i}), i=1,\hdots,n, j=1,\hdots,M,\\
(B_{A_1})_{ij} &= \lf.\calL \|\vx - \vx_{\calR^1_j} \|^m \rt|_{\vx = \vx_{\calI^1_i}}, i=1,\hdots,n, j=1,\hdots,p_1, \\
(B_{\Psi_1})_{ij} &= \lf.\calL \psi^1_i(\vx)\rt|_{\vx = \vx_{\calR^1_j}}, i=1,\hdots,M, j=1,\hdots,p_1,\\
O_{ij} &= 0, i,j = 1,\hdots,M.
\end{align}
$W_1$ is the \emph{matrix} of overlapped RBF-FD weights, with each column containing the RBF-FD weights for a point $\vx \in \mathbb{B}_1$. The linear system \eqref{eq:rbf_linsys} has a unique solution if the nodes in $P_1$ are distinct~\cite{Fasshauer:2007,Wendland:2004}. The matrix of polynomial coefficients $\Lambda_1$ is merely a set of Lagrange multipliers that enforces the polynomial reproduction constraint \eqref{eq:poly_constraint}.

\revone{
\subsection{Algorithmic description of differentiation matrix assembly}
\label{sec:assembly}
\begin{algorithm}
\caption{Differentiation matrix assembly using Overlapped RBF-FD}
\label{alg:assembly}
\begin{algorithmic}[1]	
  \Statex{\bf Given}: $X = \{\vx_k\}_{k=1}^N$, the set of nodes in the domain.
	\Statex{\bf Given}: $\delta \in (0,1]$, the overlap parameter.
	\Statex{\bf Given}: $\calL$, the linear differential operator to be approximated.
	\Statex{\bf Given}: $n << N$, the stencil size.
	\Statex{\bf Generate}: $L$, the $N \times N$ differentiation matrix approximating $\calL$ on the set $X$.
	\Statex{\bf Generate}: $N_{\delta}$, the number of stencils.
	\State Build a kd-tree on the set $X$ in $O(N \log N)$ operations.
	\State Initialize $g$, an $N$-long array of flags, to 0.
	\State Initialize the stencil counter, $N_{\delta}=0$.
	\For {$k=1,N$}
		\If {g(k) == 0}		
				\State Use kd-tree to get $\{\vx_{\calI^k_1},\hdots,\vx_{\calI^k_n}\}$. Here, $\calI^k_1 = k$.
				\State Get $R_k = \{\calR^k_1, \hdots, \calR^k_{p_k}\}$ and $\mathbb{B}_k = \{\vx_{\calR^k_1}, \hdots, \vx_{\calR^k_{p_k}} \}$ using \eqref{eq:stencil_rad}--\eqref{eq:ball_pts}.
				\State Use \eqref{eq:rbf_linsys} to compute $W_k$, the $n \times p_k$ matrix of RBF-FD weights.
				\For {$i=1,p_k$}
					\State Set $g\lf(\calR^k_i\rt) = 1$.	
					\For{$j=1,n$}
							\State Set $L\lf(\calR^k_i,\calI^k_j\rt) = W_k(j,i)$.
					\EndFor
				\EndFor
				\State Set $N_{\delta} = N_{\delta}+1$.
		\EndIf
	\EndFor
\end{algorithmic}
\end{algorithm}
We next describe how to assemble an $N \times N$ sparse differentiation matrix from the mathematical description above. The algorithm described here is a \emph{serial} one; parallel algorithms for assembly are a subject of ongoing research. The procedure involves using a kd-tree for nearest neighbor searches, together with a form of breadth-first search, and is described in Algorithm \ref{alg:assembly}. For simplicity, algorithm \ref{alg:assembly} describes the assembly process as directly populating the entries of $L$. In practice, one stores the non-zero elements, and row and column indices/pointers (depending on the sparse matrix format in use), and uses these arrays in place of $L$.

It is also worth briefly examining the cost of this process. Since Algorithm \ref{alg:assembly} uses a kd-tree, it has a preprocessing cost of $O(N \log N)$. Each iteration in the outer loop of the algorithm (line $4$) incurs the following costs: $O\lf((n-1) \log N\rt)$ to find the nearest neighbors, and $O\lf( (n+M)^3 + p_k (n+M)^2 \rt)$ for finding the weight matrix $W_k$. There are $N_{\delta}$ such steps. If the points are quasi-uniform, then $N_{\delta} \approx \frac{N}{p}$, where $p = (1-\delta)^d n$~\cite{ShankarJCP2017}. If $\delta <1$, each iteration of this loop is more expensive than the standard RBF-FD method, but there are far fewer iterations (by a factor $\frac{1}{p}$). The assembly process as described above also requires fewer nearest neighbor searches than the standard RBF-FD method if $\delta <1$.
}

\subsection{Local error estimates}
\label{sec:error_estimate}
The polynomial reproduction property \eqref{eq:poly_constraint} of the augmented local RBF interpolant can be used to develop a local error estimate~\cite{DavydovMinimal2016}. Following~\cite{ShankarJCP2017}, if $\calL$ is a differential operator of order $\theta$ applied to some function $f$, this error estimate can be written as:
{\small
\begin{align}
\lf| \calL f(\vy) - \sum\limits_{j=1}^n w_j(\vy) f(\vx_j) \rt| \leq P(\vy)\max\limits_{|\alpha| = \ell} \|D^{\alpha} f\|_{W^{1,\infty}(\Omega_1)} \lf( h(\vy) \rt)^{\ell+1-\theta}, 
\label{eq:error_estimate}
\end{align}
}
where $\vy$ is some evaluation point, $h(\vy)$ is the largest distance between the point $\vy$ and every point $\vx$ in the stencils, $P(\vy)$ is some \emph{growth function}~\cite{DavydovMinimal2016}, and $\alpha$ is a multiindex. The term $\max\limits_{|\alpha| = \ell} \|D^{\alpha} f\|_{W^{1,\infty}(\Omega_1)}$ is simply the Sobolev $\infty$-norm of $f$ on the stencil $P_1$ with convex hull $\Omega_1$. The error estimate shows \revone{shows that the order of convergence in any RBF-FD method depends primarily on the degree $\ell$, not the location of the evaluation point $\vy$. One can therefore expect the same orders of convergence for the overlapped RBF-FD method as in the standard one, assuming $\vy$ is not too close to the boundary of $\Omega_1$. As was shown in~\cite{ShankarJCP2017}, the difference in accuracy between the standard and overlapped RBF-FD methods is virtually non-existent for higher-order approximations, can be kept small for lower-order approximations using judicious choices of $\delta$.}

\subsection{Parameter selection}
We now describe how the different parameters in the overlapped RBF-FD method are selected.  Based on \eqref{eq:error_estimate}, we know the polynomial degree $\ell$ controls the order of the approximation for a differential operator of order $\theta$. Thus, given a differential operator of order $\theta$, if we require an RBF-FD method with order of accuracy $\xi$, we set
\begin{align}
\ell = \xi + \theta - 1.
\label{eq:select_ell}
\end{align}
For instance, if we are approximating the gradient operator $\nabla$ ($\theta = 1$) to third-order accuracy ($\xi = 3$), this gives $\ell = 3$. On the other hand, if we are approximating the Laplacian operator $\Delta$ ($\theta = 2$) to third-order accuracy ($\xi = 3$), this gives $\ell = 4$. The stencil size $n$ in $d$ dimensions can be deduced from the relationships $M = {\ell + d \choose d}$ and $M \leq \lf \lfloor\frac{n}{2} \rt\rfloor$ to be $n \geq 2M$. In practice, it appears to be beneficial to use a few more nodes per stencil than $2M$, though the optimal choice is likely problem dependent. We use
\begin{align}
n = 2M + \lf \lfloor \ln(2M) \rt\rfloor.
\label{eq:select_n}
\end{align}
We must also choose the overlap parameter $\delta \in (0,1]$. \revone{For small values of $n$ (equivalently of $\ell$), this choice trades off speedup with accuracy as $\delta \to 0$; however, as $n$ increases, smaller values of $\delta>0.2$ barely impact accuracy~\cite{ShankarJCP2017}}. In practice, we have observed that setting $\delta \leq 0.2$ typically completely decouples the stencils, resulting in ill-posed subproblems when solving PDEs. Given these constraints, we use the following heuristic:
\[
 \delta =
  \begin{cases} 
      \hfill 0.7 \hfill & \text{ if $\ell \leq 3$} \\
      \hfill 0.5 \hfill & \text{ if $4 \leq \ell < 6$} \\
			\hfill 0.3 \hfill & \text{ if $\ell \geq 6$} \\
  \end{cases}
\]
The automatic selection of $\delta$ based on stability properties of augmented local RBF interpolants is a subject of ongoing research. It is likely that the optimal choice of $\delta$ depends on the local Lebesgue functions corresponding to the differential operator $\calL$; see Section \ref{sec:stab} and~\cite{ShankarJCP2017,DavydovMinimal2016} for discussions on these local Lebesgue functions.

\revone{
\textbf{Connection to other RBF methods}
As described here, the overlapped RBF-FD method bears some resemblance to both the standard RBF-FD method and the RBF-Partition of Unity (RBF-PU) method~\cite{Wendland:2004,Fasshauer:2007}. In the latter case, one computes a \emph{global} interpolant to the node set $X$ by dividing $X$ into a set of overlapping partitions called \emph{patches} (controlled by an overlap parameter), forming a local RBF approximation on each patch, and then blending these approximations together with weight functions that are compactly-supported on each patch. As such, the standard RBF-FD method can be thought of as a very special RBF-PU method with weight functions that are exactly 1 at each node, and zero elsewhere, resulting in as many patches as there are nodes. The overlapped RBF-FD method can then be thought of as an RBF-PU method where the weight functions are 1 on the balls $\mathbb{B}_k$, zero elsewhere, resulting in fewer patches than nodes. However, in our opinions, the overlapped RBF-FD method is best viewed as an RBF-FD method where one makes a choice of either centered or one-sided approximations for each node in $X$, thereby resulting in stencil-sharing for some nodes. Indeed, the overlapped and standard RBF-FD methods share the same local error estimates from Section \ref{sec:error_estimate}. That said, it may be that viewing the assembly process from the RBF-PU perspective may help with parallelizing our method. This is a topic of ongoing research.
}
\section{A new scaling law for $m$ and $\ell$}
\label{sec:stab}

\begin{figure}[h!]
\centering
\subfloat[$m=2\ell+1$]
{
	\includegraphics[scale=0.6]{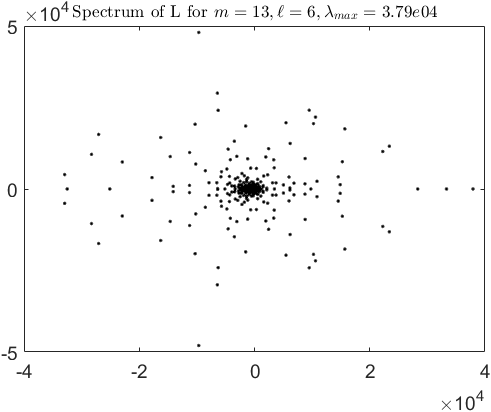}
	\label{fig:c11}	
}
\subfloat[$m=\ell+1$]
{
	\includegraphics[scale=0.6]{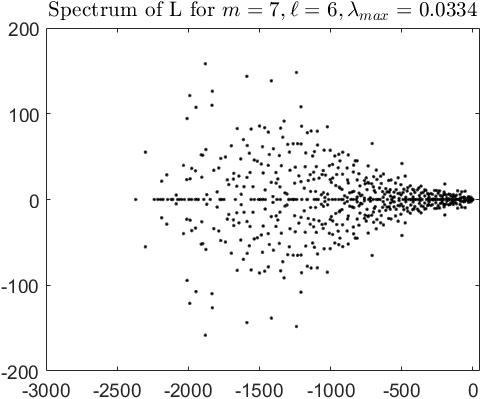}
	\label{fig:c12}	
}

\subfloat[$m = \ell-1$]
{
	\includegraphics[scale=0.6]{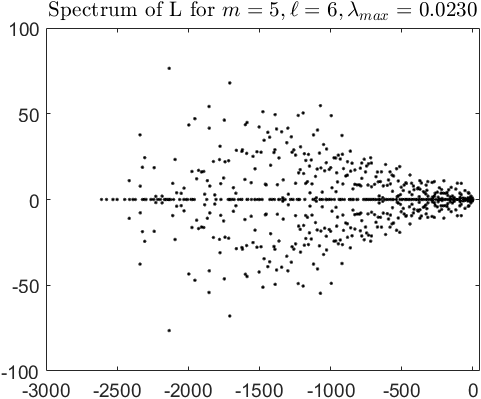}
	\label{fig:c21}	
}
\subfloat[$m = 3$]
{
	\includegraphics[scale=0.6]{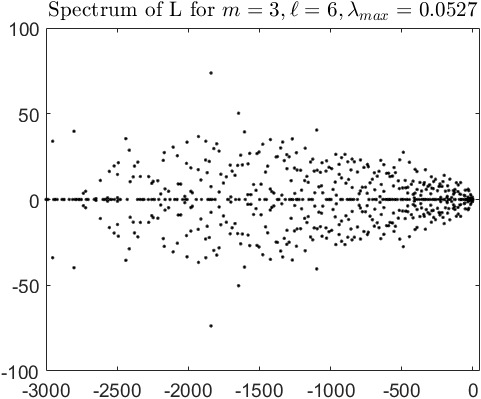}
	\label{fig:c22}	
}
\caption{Spectrum of the discrete Laplacian $L$ on a node set with $N=658$ nodes for different relationships between $m$ and $\ell$. The node set here was a set of Poisson disk samples generated on the disk using the method from~\cite{SFKSISC2017}.}
\label{fig:eigs_laws1}
\end{figure}

The selection of the degree $m$ of the PHS RBF requires more thought. On the one hand, $m$ is known to control the constant in the RBF-FD error estimate \eqref{eq:error_estimate}, with higher values of $m$ giving lower errors~\cite{FlyerPHS,BarnettPHS}. On the other hand, larger values of $m$ may cause greater instability, as the smallest eigenvalue of the block matrix from \eqref{eq:rbf_linsys} can be bounded from below by $q^m$, where $q$ is the separation distance of the stencil $P_1$; see~\cite{Wendland:2004} for some preliminary estimates, and~\cite{FlyerNS} for remarks about stability. A third more subtle issue is the question of whether $m$ and $\ell$ should be related. Historically, it was common to use the \emph{classical scaling law}~\cite{iske2002}:
\begin{align}
m = 2\ell+1,
\label{eq:scaling_law_classical}
\end{align}
where $\ell$ is selected to be the minimum to prove unisolvency for the RBF interpolant with $\phi(r) = r^m$. This approach was also recently used for PHS-based methods for advection on the sphere~\cite{ShankarWrightJCP2017}, and advection and reaction-diffusion equations on more general manifolds~\cite{SNKJCP2018}. However, modern PHS RBF-FD methods for Euclidean domains decouple $m$ and $\ell$~\cite{FlyerNS,FlyerPHS,BarnettPHS,ShankarJCP2017}. Our goal is to attempt to find a scaling law that produces the smallest spurious eigenvalue (if any) in the spectrum of the discrete Laplacian $L$. To that end, we explored two alternatives to the classical scaling law. The first alternative is:
\begin{align}
m = 
\begin{cases}
	\hfill \ell \hfill & \text{if $\ell$ is odd} \\
	\hfill \ell+1 \hfill & \text{if $\ell$ is even}\\
\end{cases}
\label{eq:scaling_law_1}
\end{align}
which serves to gently increase $m$ as $\ell$ is increased. The second alternative, which we ultimately settled on, is:
\begin{align}
m = 
\begin{cases}
	\hfill \ell \hfill & \text{if $\ell$ is odd} \\
	\hfill \ell-1 \hfill & \text{if $\ell$ is even}\\
\end{cases}
\label{eq:scaling_law_2}
\end{align}
The eigenvalues for the discrete Laplacian $L$ computed using the classical scaling law, these two alternatives, and the choice of $m=3$ are shown for $N=658$ nodes on the disk using $\ell=6$ in Figure \ref{fig:eigs_laws1}.
From Figure \ref{fig:c11}, it is clear that the classical scaling law $m=2\ell+1$ results in a large spurious eigenvalue. In general, all the other scaling laws appear to give very similar results, albeit with some subtle differences. Figure \ref{fig:c12} uses the scaling law $m=\ell+1$, while Figure \ref{fig:c21} uses the scaling law $m = \ell-1$; clearly, both produce almost identical results, with the exception that $m=\ell+1$ results in a slightly larger spurious eigenvalue. Continuing this trend, it may be tempting to always simply choose a very small reasonable value of $m$. However, Figure \ref{fig:c22} with $m=3$ shows that this may not be entirely advisable either, since the largest spurious eigenvalue here is actually larger than the cases of $m=\ell+1$ and $m=\ell-1$. 
\begin{figure}[h!]
\centering
\subfloat[$m=2\ell+1$]
{
	\includegraphics[scale=0.6]{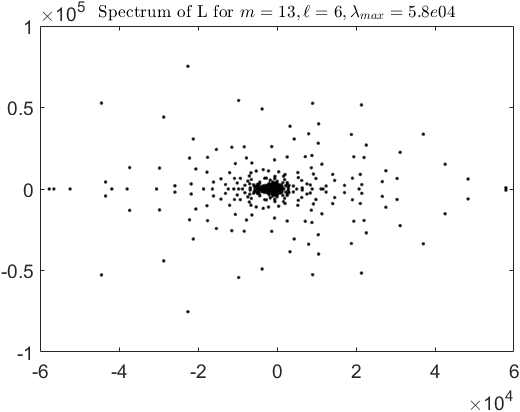}
	\label{fig:d11}	
}
\subfloat[$m=\ell+1$]
{
	\includegraphics[scale=0.6]{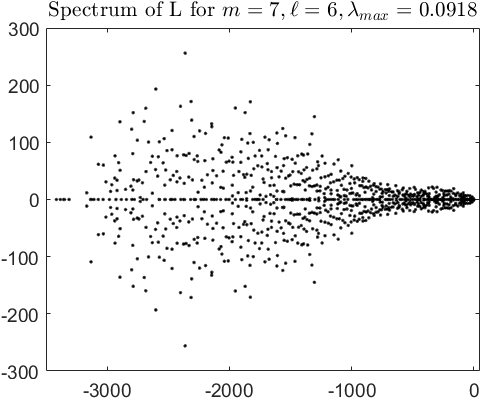}
	\label{fig:d12}	
}

\subfloat[$m = \ell-1$]
{
	\includegraphics[scale=0.6]{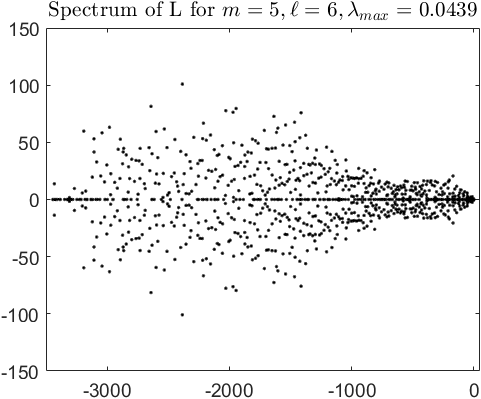}
	\label{fig:d21}	
}
\subfloat[$m = 3$]
{
	\includegraphics[scale=0.6]{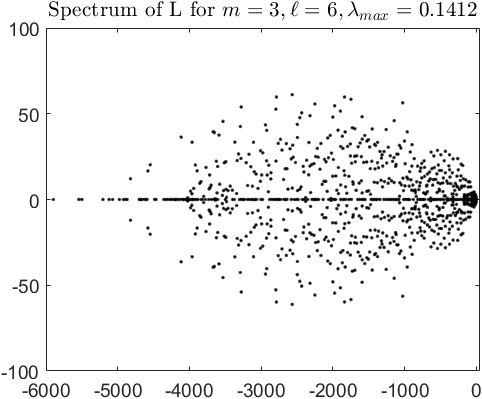}
	\label{fig:d22}	
}
\caption{Spectrum of the discrete Laplacian $L$ on a node set with $N=919$ nodes for different relationships between $m$ and $\ell$. The node set here was a set of Poisson disk samples generated on the disk using the method from~\cite{SFKSISC2017}.}
\label{fig:eigs_laws2}
\end{figure}
The spectrum of $L$ for different scaling laws for $N=919$ nodes with $\ell=6$ is shown in Figure \ref{fig:eigs_laws2}. Once again, $m=2\ell+1$ can be rejected as extremely unstable, even after refinement (Figure \ref{fig:d11}). However, the gaps between the other scaling laws have widened slightly as well. For instance, $m=\ell+1$ (Figure \ref{fig:d12}) now results in a spurious eigenvalue that is twice as large as $m=\ell-1$ (Figure \ref{fig:d21}). Even more interestingly, $m=3$ (Figure \ref{fig:d22}) results in a spurious eigenvalue whose real part is much larger than the $m=\ell+1$ and $m=\ell-1$ cases. Once again, the most stable choice appears to be $m = \ell-1$. Picking a fixed value of $m=3$ (or any other number) regardless of the choice of $\ell$ may be inadvisable on Euclidean domains except when the node sets are very regular. These trends continue for different values of $\ell$ (not shown). It may be the case that the correct scaling relationship needs to be reevaluated for each type of node set used.

\begin{figure}[h!]
\centering
\subfloat[]
{
	\includegraphics[scale=0.4]{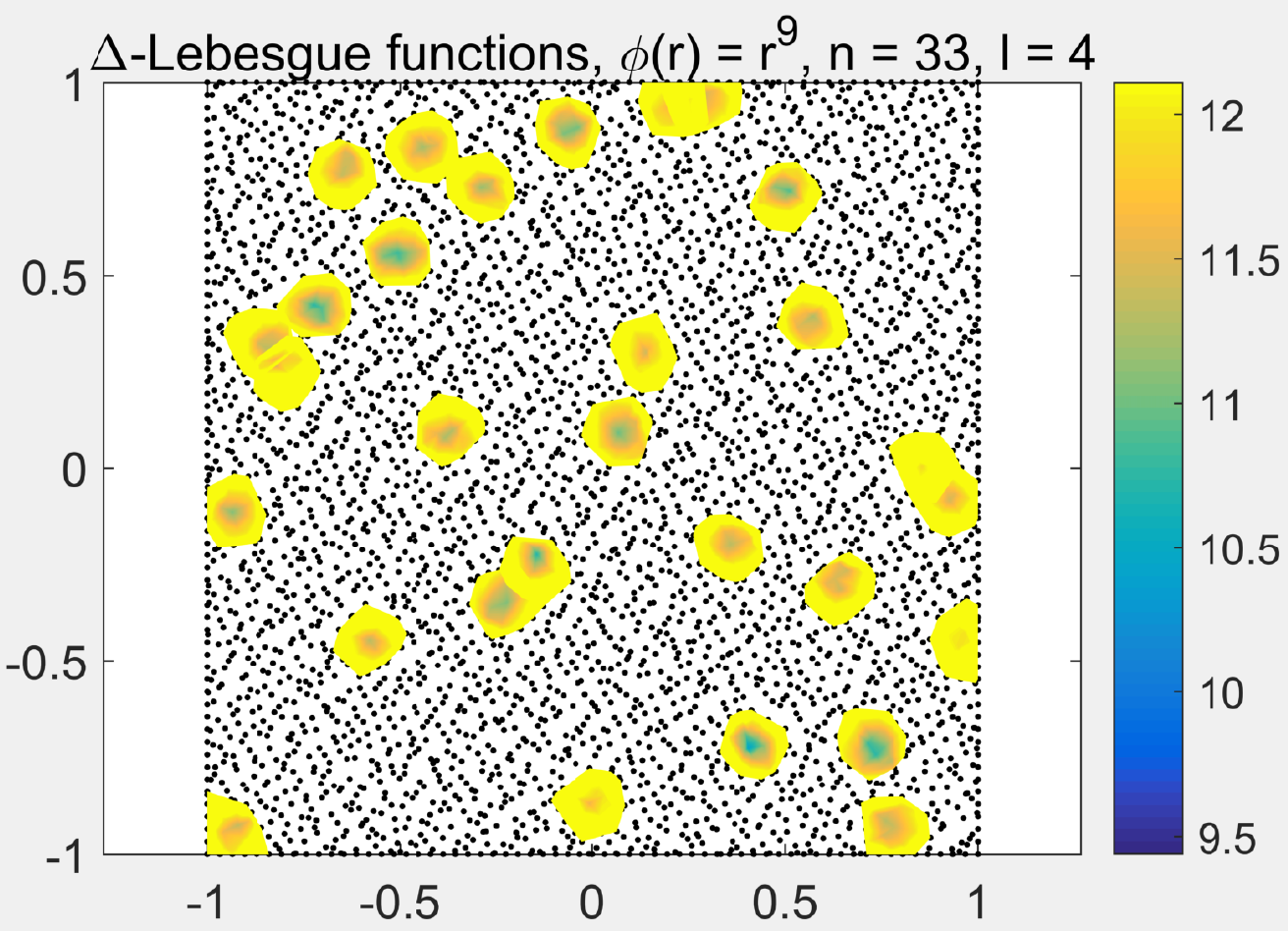}
	\label{fig:leb_c13}	
}
\subfloat[]
{
	\includegraphics[scale=0.4]{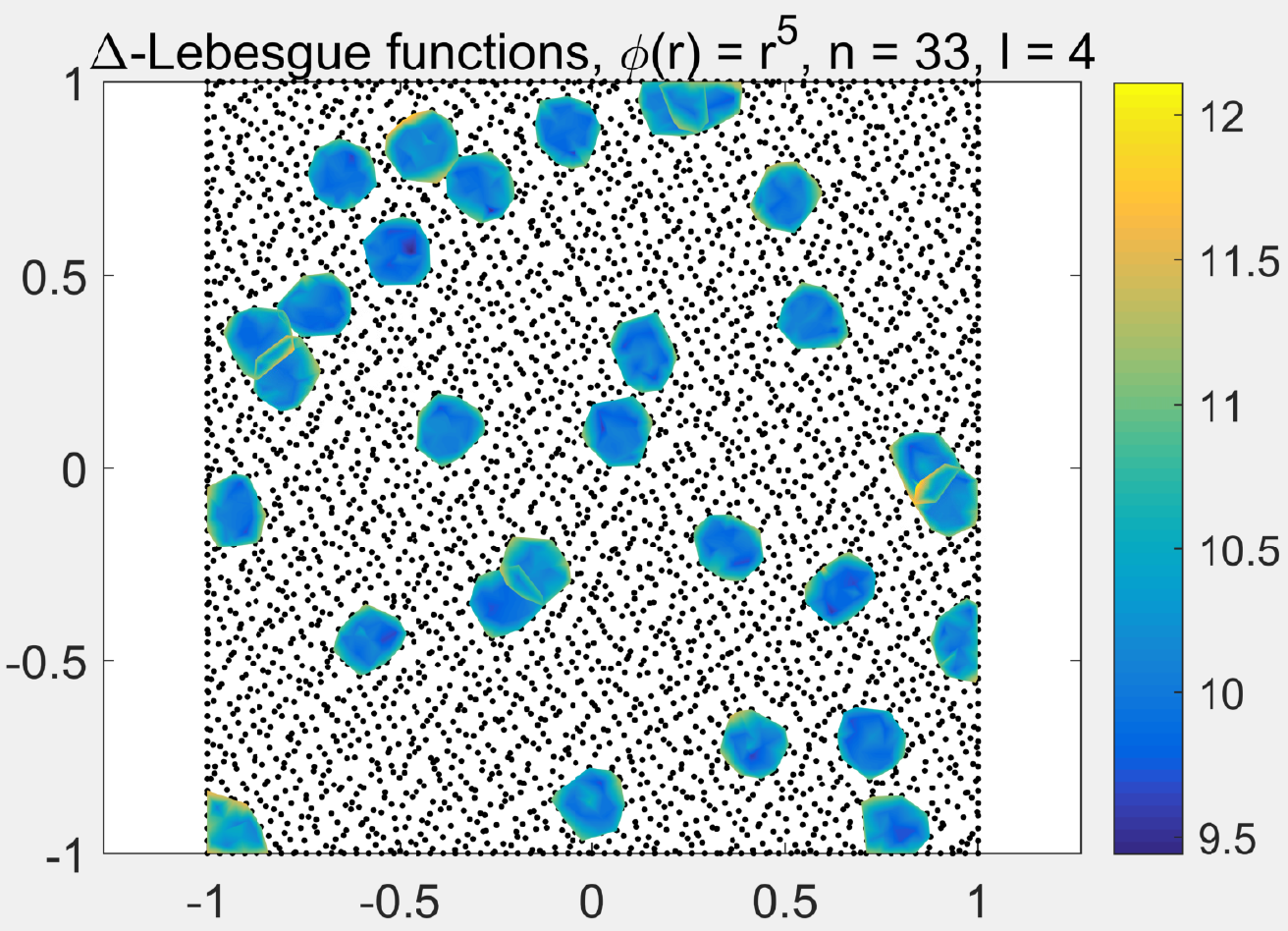}
	\label{fig:leb_c11}	
}
\subfloat[]
{
	\includegraphics[scale=0.4]{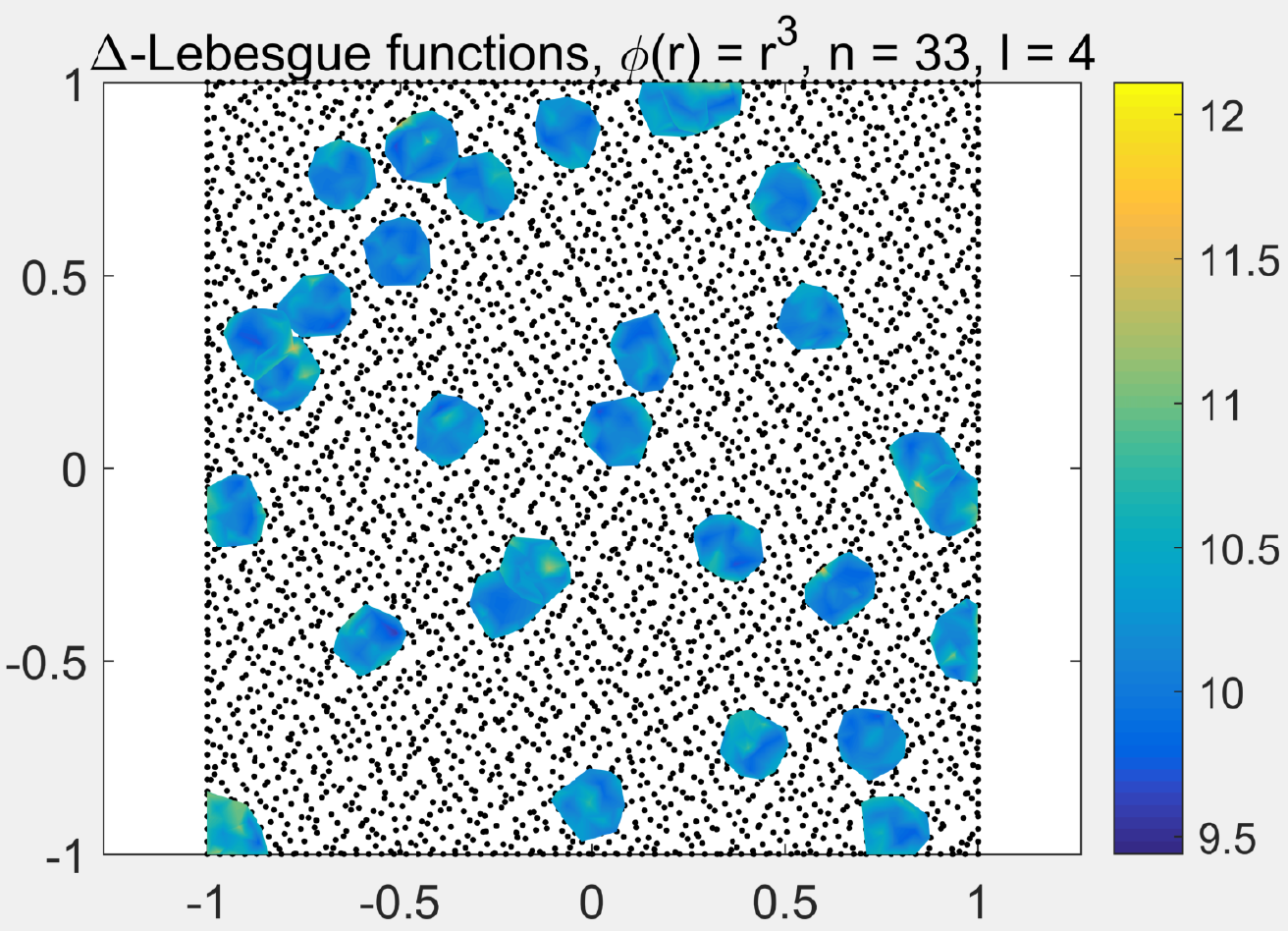}
	\label{fig:leb_c12}	
}

\subfloat[]
{
	\includegraphics[scale=0.4]{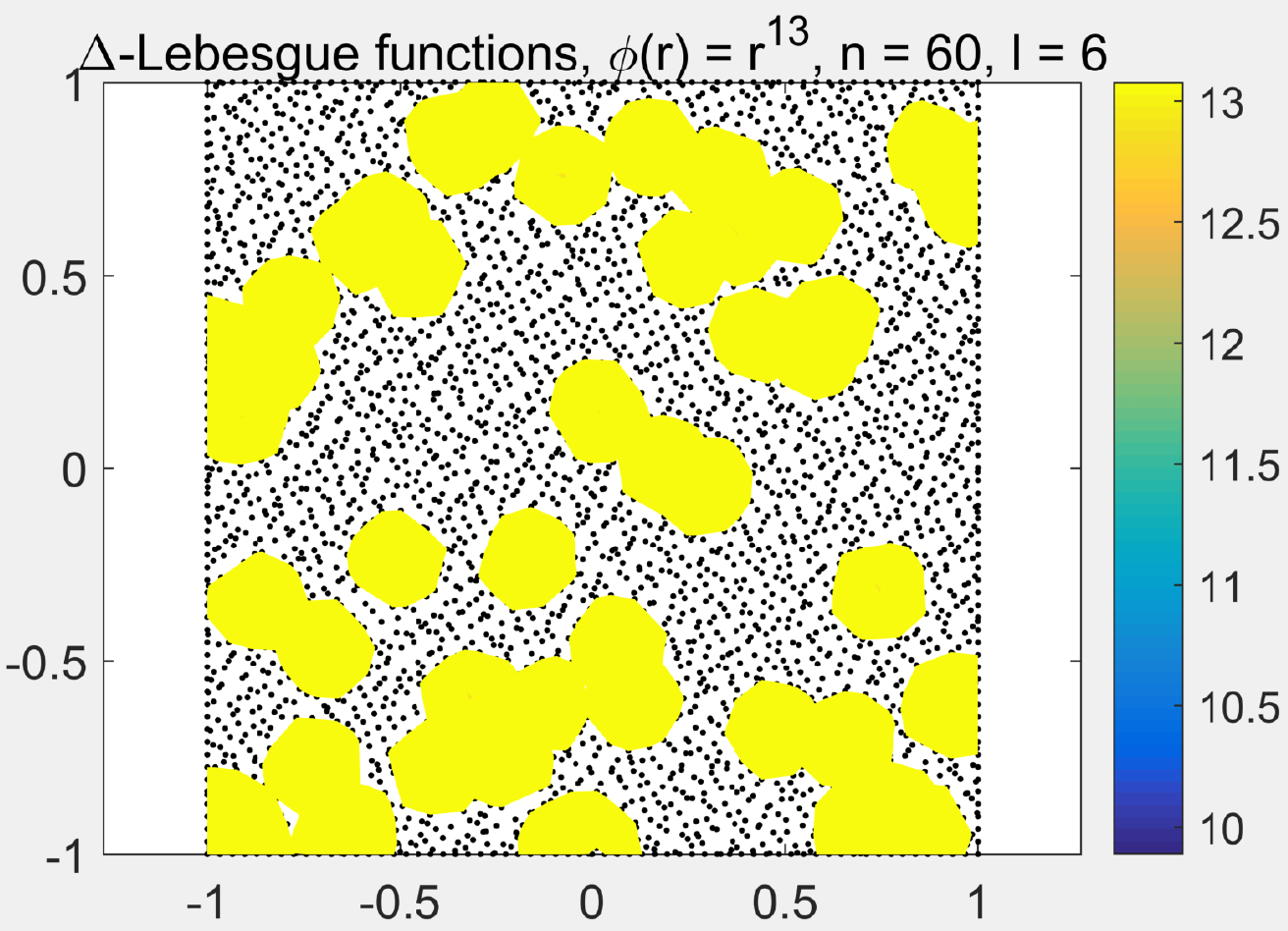}
	\label{fig:leb_c23}	
}
\subfloat[]
{
	\includegraphics[scale=0.4]{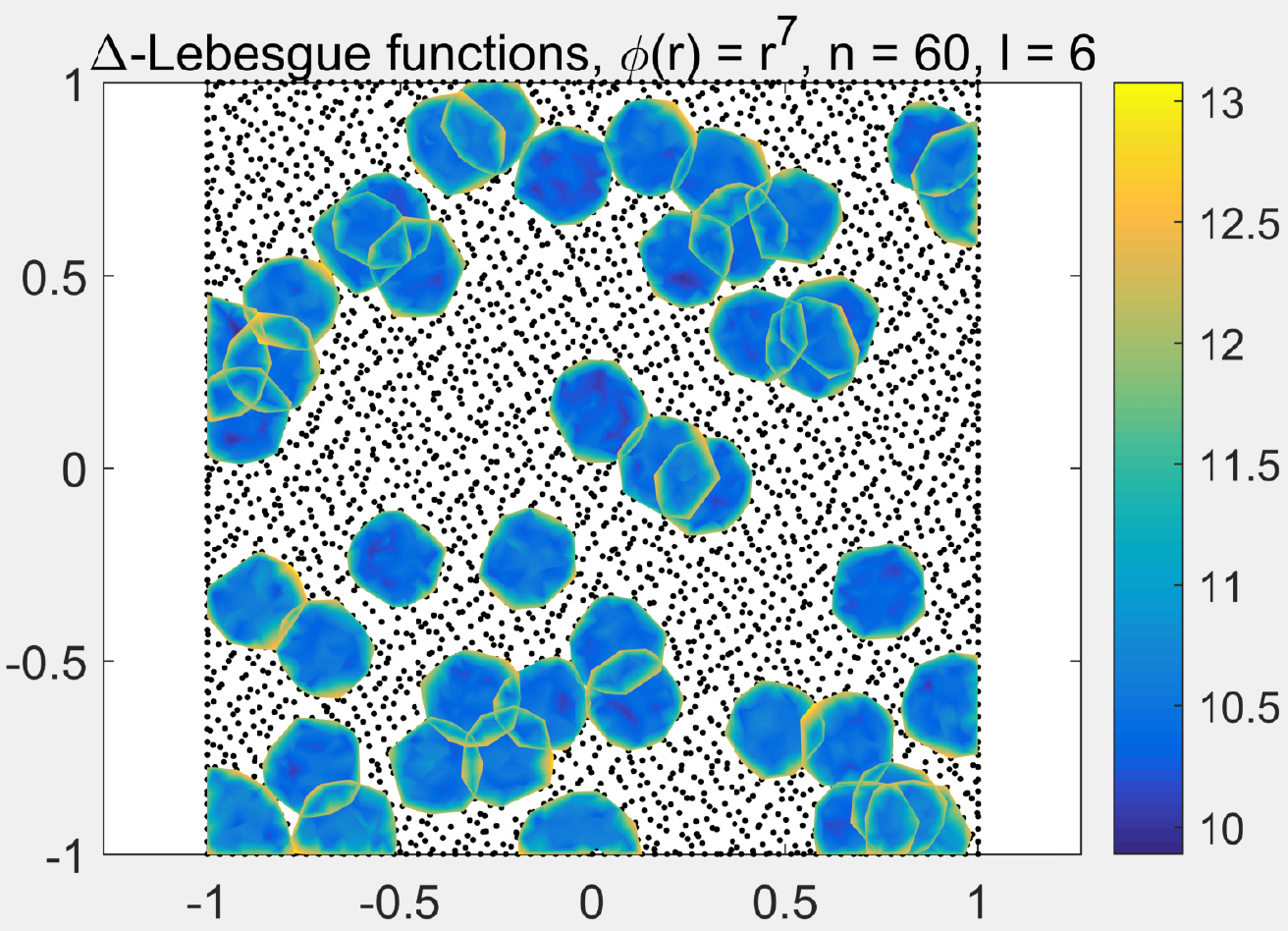}
	\label{fig:leb_c21}	
}
\subfloat[]
{
	\includegraphics[scale=0.4]{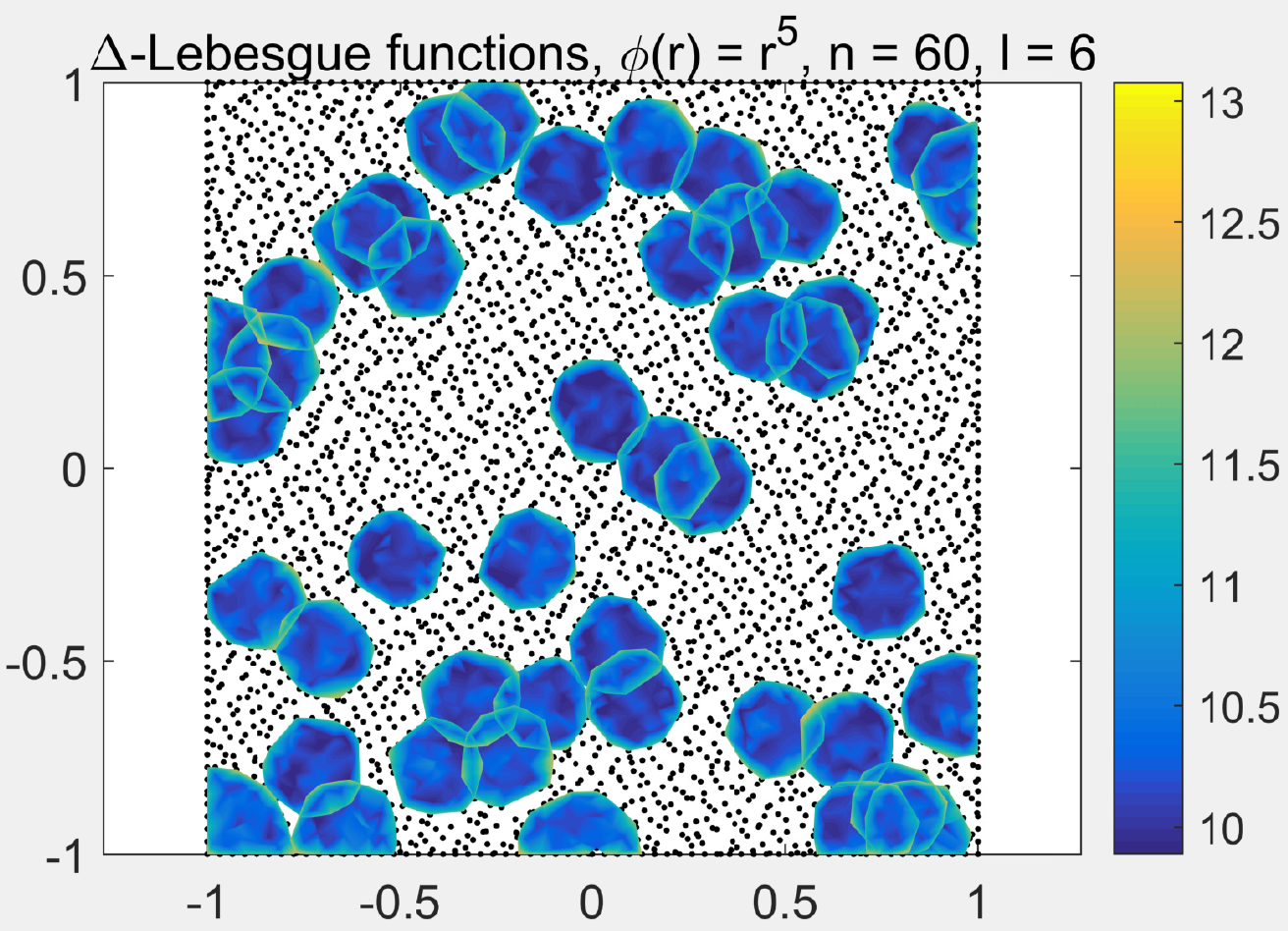}
	\label{fig:leb_c22}	
}

\subfloat[]
{
	\includegraphics[scale=0.4]{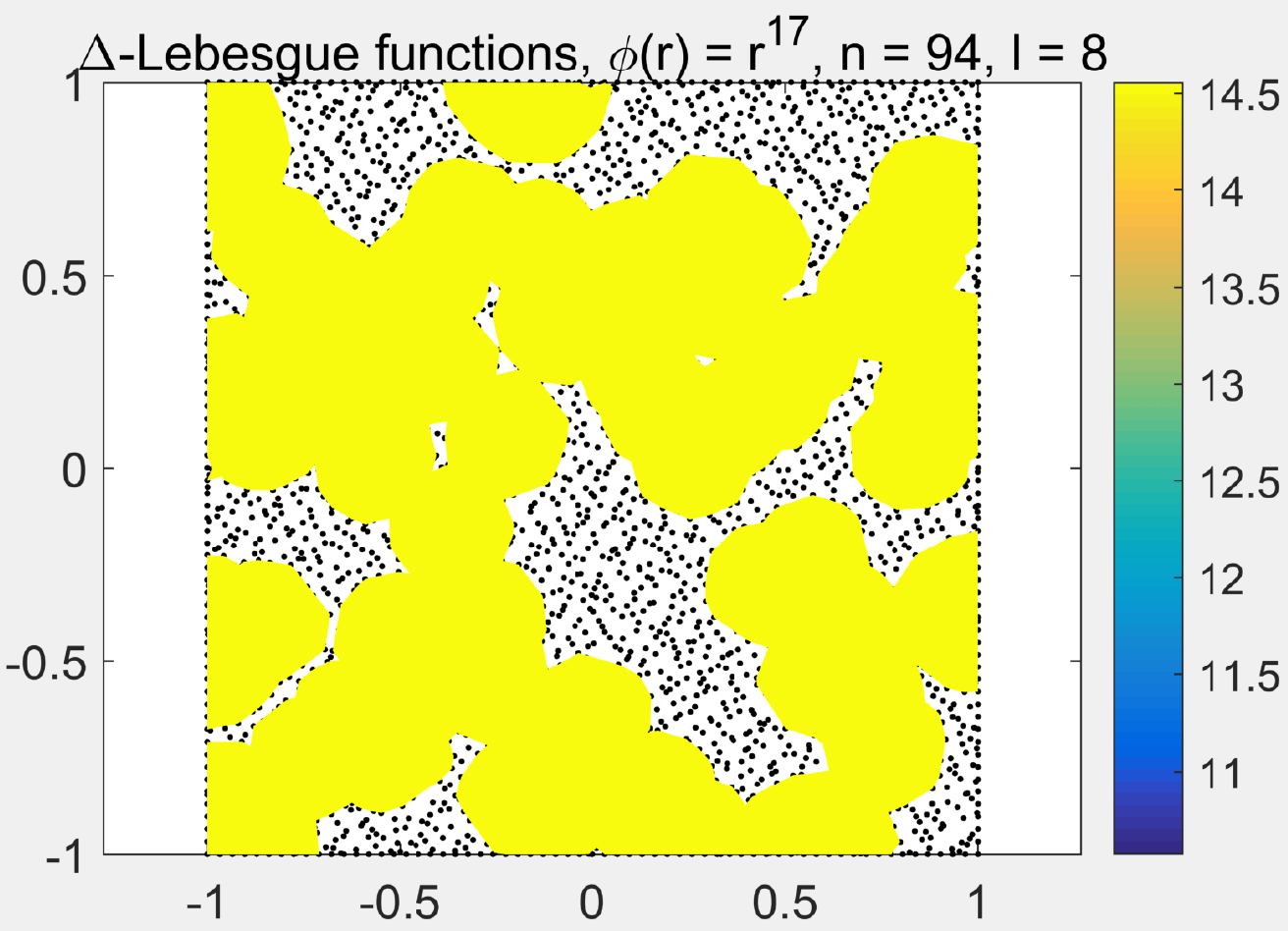}
	\label{fig:leb_c33}	
}
\subfloat[]
{
	\includegraphics[scale=0.4]{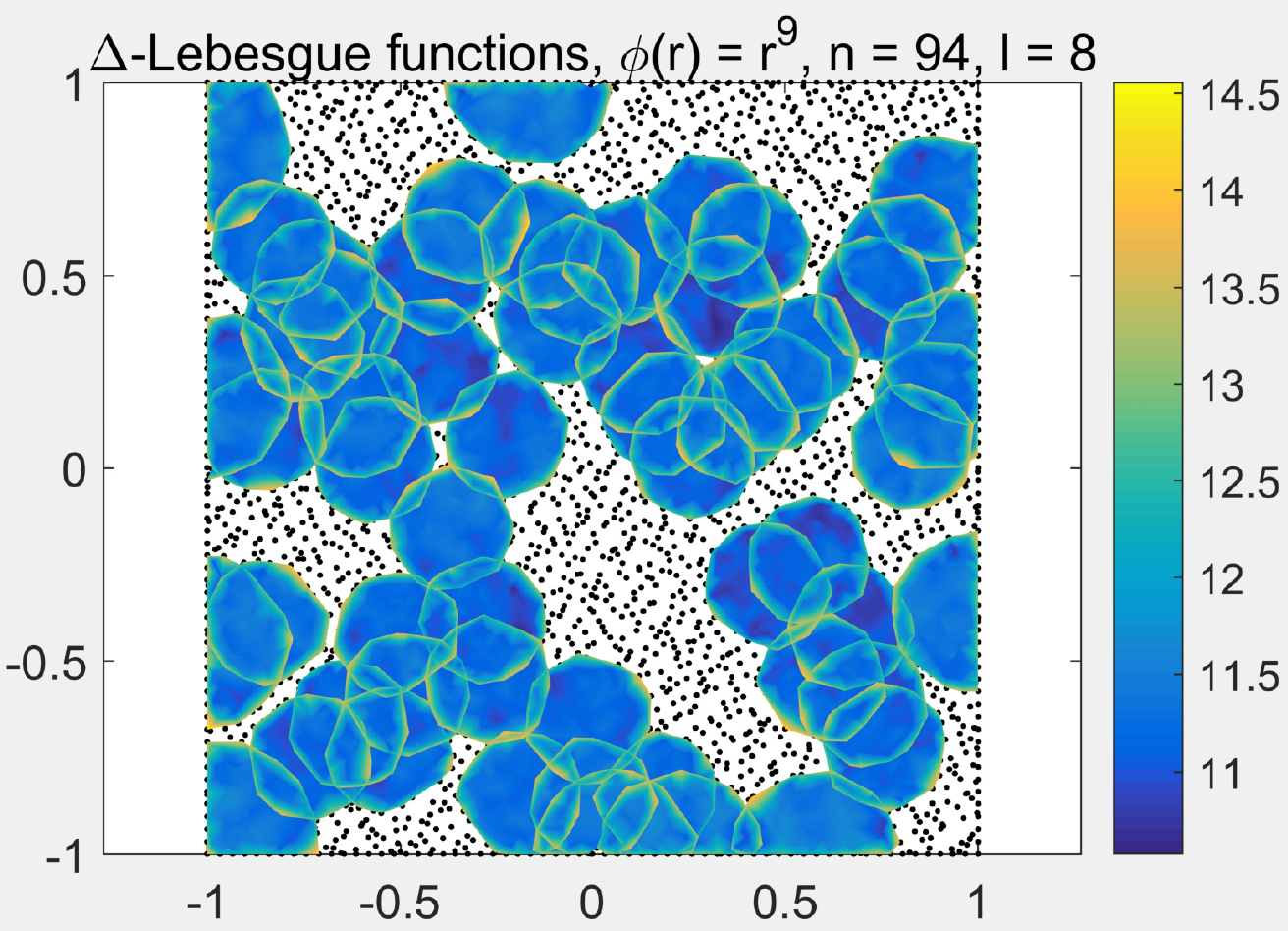}
	\label{fig:leb_c31}	
}
\subfloat[]
{
	\includegraphics[scale=0.4]{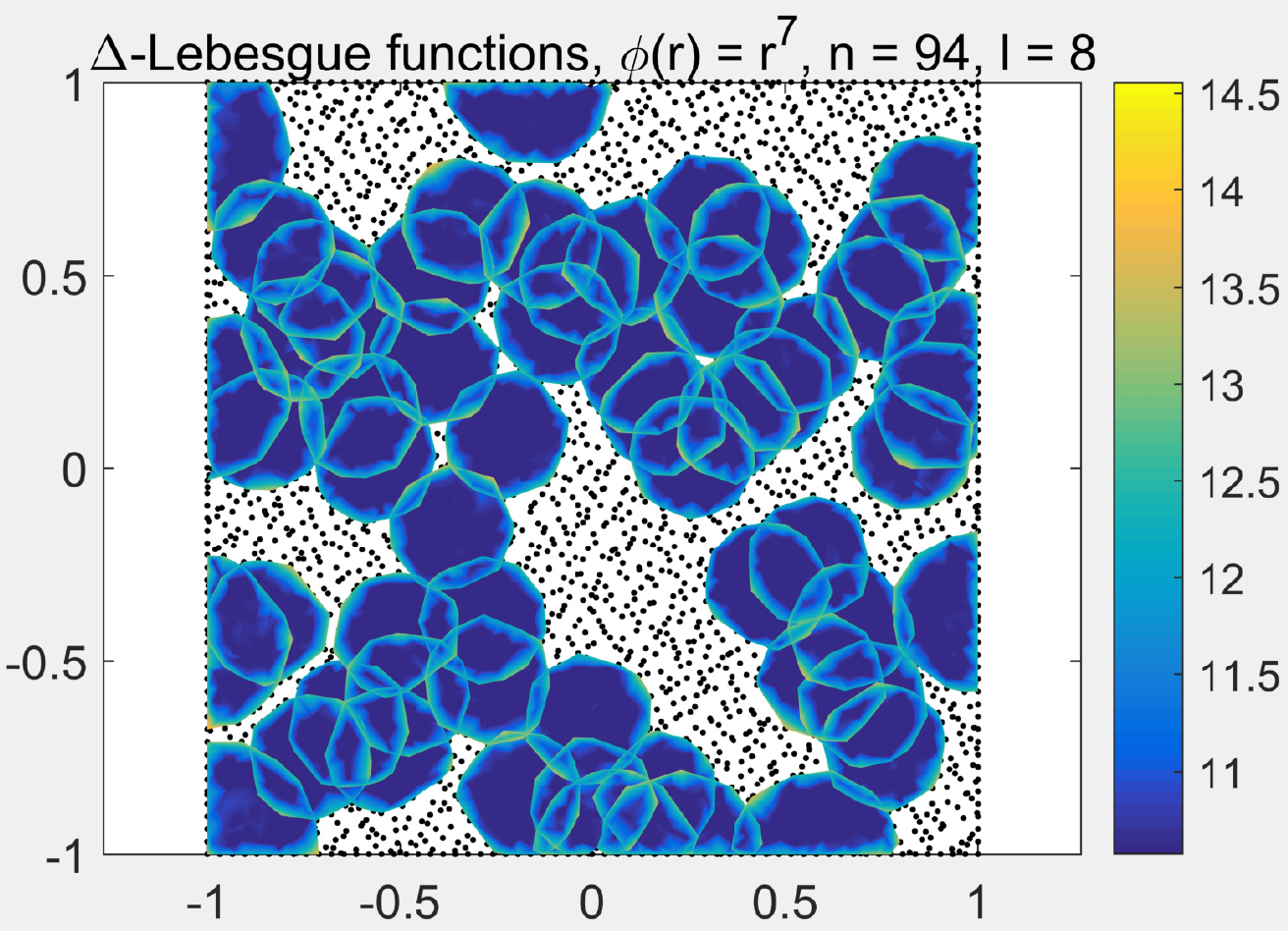}
	\label{fig:leb_c32}	
}
\caption{Local $\Delta$-Lebesgue functions for different scaling laws on a logarthmic scale. The figure contrasts \eqref{eq:scaling_law_classical} (left), \eqref{eq:scaling_law_1} (middle), and \eqref{eq:scaling_law_2} (right). Lighter colors indicate larger values, which in turn indicate greater potential for instability.}
\label{fig:leb_c1}
\end{figure}

Thus far, we have only discussed the \emph{global} stability associated with the different scaling laws. It is reasonable to expect that these global patterns are reflected in some way at the stencil level. In other words, we expect that the global stability must be connected to \emph{local} stability in some way.  Before we can explore this connection, we must first define a few quantities. Returning to the stencil $P_1$, consider the \emph{local $\calL$-Lebesgue function} of the augmented RBF interpolant $s_1(\vx,\vy)$ used for approximating the operator $\calL$, defined as
\begin{align}
\Lambda^1_{\calL}(\vy) = \sum\limits_{j=1}^n |w^1_j(\vy)| = \|{\bf w}^1\|_1,
\end{align}
that is, the 1-norm of the RBF-FD weight vector at a point. Let $w^1_1(\vy)$ be the RBF-FD weight associated with the stencil point $\vy$. From~\cite{ShankarJCP2017}, we have the sufficient condition for all eigenvalues of the differentiation matrix $L$ (corresponding to $\calL$) having non-positive real parts as:
\begin{align}
2w^1_1(\vy) \geq - \Lambda^1_{\calL}(\vy).
\end{align}
This condition highlights the connection between the local Lebesgue function (itself connected to the stability of the local RBF interpolant~\cite{Wendland:2004}) and the \emph{global} stability of the RBF-FD method. In general, large values of $\Lambda^1_{\calL}(\vy)$ are associated with instability, and typically produce spurious eigenvalues with positive real parts~\cite{ShankarJCP2017}. To examine the impact of the scaling laws \eqref{eq:scaling_law_classical}--\eqref{eq:scaling_law_2} on the local Lebesgue functions (and therefore local stability), we will now study the local $\Delta$-Lebesgue functions for a set of $N = 4000$ Halton points in $[-1,1]^2$, where $\Delta$ is the Laplacian; evenly-spaced points are used on the boundary; these Halton points are typically even more scattered than Poisson disk samples, and represent a worst-case scenario. The results are shown in Figure \ref{fig:leb_c1}.
The first column in Figure \ref{fig:leb_c1} corresponds to the classical law \eqref{eq:scaling_law_classical}, the second to \eqref{eq:scaling_law_1}, and the third to \eqref{eq:scaling_law_2}. The colorbar scale is fixed across each row, but allowed to increase between rows to make clear that RBF-FD weights can increase in magnitude as the polynomial degree is increased. Immediately, it is clear that the classical scaling law results in the largest weights, and is consequently almost always unstable (Figures \ref{fig:leb_c13},\ref{fig:leb_c23}, \ref{fig:leb_c33}).  Careful observation also shows that the last column (corresponding to \eqref{eq:scaling_law_2}) results in more stable approximations than the middle column (corresponding to \eqref{eq:scaling_law_1}), which matches our conclusions from the eigenvalue plots (Figures \ref{fig:eigs_laws1} and \ref{fig:eigs_laws2}). In general, we find that the scaling law \eqref{eq:scaling_law_2} produces stable approximations across a wide range of parameters, and therefore we use it in this article.

\section{Hyperviscosity-based Stabilization}
\label{sec:robust_hyp}

As mentioned previously, the RBF-FD differentiation matrices corresponding to the discretized gradient operator (and occasionally even the Laplacian operator) can contain eigenvalues with positive real parts. Such eigenvalues can cause spurious growth in the numerical solutions of advection-diffusion equations. Our approach to ameliorating this issue is to add a small vanishing amount of artificial hyperviscosity to the model. This transforms Equation 1 to
\begin{align}
\frac{\partial c}{\partial t} &= -\vu \cdot \nabla c + \nu \Delta c + \gamma \Delta^k c,
\label{eq:adr_hyp}
\end{align}
where $\gamma \Delta^k c$ is the artificial hyperviscosity term. This approach has been used in conjunction with explicit time-stepping in convective PDEs~\cite{FlyerNS,FoL11,BarnettPHS,FlyerLehto2012}. Our goals in this section are twofold: first, to present expressions for $\gamma$ and $k$, and second, to discuss the spatial approximation of the hyperviscosity term in the context of the RBF-FD method. We proceed by first deriving expressions for $\gamma$ in the advection-dominated regime, in the diffusion-dominated regime, then for the case where the two processes co-exist.

\subsection{Selecting $\gamma$ in the advection-dominated regime ($\nu \to 0$)}
\label{sec:gamma_1}

As mentioned in Section \ref{sec:intro}, the literature on modern RBF-FD methods gives a prescription for $\gamma$. In~\cite{FlyerNS,BarnettPHS}, for instance, Flyer et al. set this parameter to $\gamma = (-1)^{1-k} 2^{-6} h^{2k-1}$, where $k$ is the power of the Laplacian. We will show with our analysis that this formula is a special case of a much more general formula.

Assume for now a simplified 1D version of \eqref{eq:ad} where $\nu = 0$. This reduces it to the advection equation
\begin{align}
\frac{\partial c}{\partial t} + u\frac{\partial c}{\partial x} = 0.
\end{align}
Assume that $c$ is a plane wave of the form $c(x,t) = \hat{c}(t) e^{i\hat{k}x}$, where $\hat{k}$ is the wave number. The above equation should simply translate the plane wave to the right at a speed $u$. Unfortunately, the problem lies in the approximation of $\frac{\partial}{\partial x}$. In general, an RBF-FD discretization will instead result in our solving the \emph{auxiliary} equation
\begin{align}
\frac{\partial c}{\partial t} + u \frac{\tilde{\partial c}}{\partial x} = 0,
\end{align}
where $\frac{\tilde{\partial c}}{\partial x}$ is an \emph{auxiliary} differential operator that introduces a new term $\ep \hat{k}^q$ into the right hand side of the above equation so that we have
\begin{align}
\frac{\partial c}{\partial t} + u i \hat{k} \hat{c}(t) e^{i\hat{k}x} = u \ep \hat{k}^q \hat{c}(t) e^{i\hat{k}x},
\end{align}
where $\ep \geq 0$ causes \emph{growth} in $c$, and $q \in \mathbb{R}$. In general, $\ep$ and $q$ are unknown \emph{a priori}, but are functions of the node set and basis functions used; see~\cite{Platte2006} for a discussion on the relationships between eigenvalues of differential operators, RBFs, and node sets. We assume in this work that an upper-bound for $\ep$ is the real part of the spurious eigenvalue with the largest real part. In the RBF-FD context, we wish to add some amount of artificial hyperviscosity to stabilize the auxiliary PDE, giving us
\begin{align}
\frac{\partial c}{\partial t} + u \frac{\tilde{\partial} c}{\partial x}=  \gamma \frac{\partial^{2k} c}{\partial x^{2k}},k \in \mathbb{N}.
\label{eq:a_hyp_1D_aux}
\end{align}
We will now derive expressions for $\gamma$ that \emph{cancel} out the spurious growth mode $\ep \hat{k}^q$.

\subsubsection{Time-stepping with explicit RK methods}
\label{sec:hyp_exp}

Our goal in this section is to derive a hyperviscosity formulation for explicit RK methods. We first examine the forward Euler discretization as this is instructive. Discretizing \eqref{eq:a_hyp_1D_aux} using the Forward Euler scheme, we obtain
\begin{align}
c^{n+1} - c^{n} + \Delta t u\frac{\tilde{\partial} c^n}{\partial x} = \Delta t\gamma \frac{\partial^{2k} c^n}{\partial x^{2k}},
\end{align}
where the superscripts indicate time levels. Now, akin to Von Neumann analysis, assume that $c^{n+1} = Gc^n$, where $G$ is some growth factor. For time stability, we require that $|G|\leq 1$. Substituting the plane wave definition of $c^n$, we obtain
\begin{align}
G \hat{c}(t_n) e^{i\hat{k}x} = \hat{c}(t_n) e^{i\hat{k}x} -\Delta t u \hat{c}(t_n)e^{i\hat{k}x} \lf(i\hat{k} -  \ep \hat{k}^q\rt) + \Delta t \gamma i^{2k} \hat{k}^{2k} \hat{c}(t_n)e^{i\hat{k}x}.
\end{align}
This simplifies to
\begin{align}
G&= \lf(1 - \Delta t u i\hat{k} \rt) + \lf(u\Delta t \ep \hat{k}^q + \Delta t \gamma i^{2k} \hat{k}^{2k} \rt).
\label{eq:paren}
\end{align}
Even in the absence of the spurious growth term $u\Delta t \ep \hat{k}^q $, we have $|G|> 1$ due to the first parenthetical term, implying that forward Euler is always unstable on the semi-discrete advection equation. This is because its stability region does not include the imaginary axis. This motivates the use of other Runge-Kutta (RK) methods that contain the imaginary axis. In general, $G(z)$ is a polynomial for all explicit RK methods. For example, for forward Euler, we have
\begin{align}
G(z) & = 1 + z,
\end{align}
while for the classical fourth-order RK method (RK4) we have
\begin{align}
G(z) & = 1 + z + \frac{1}{2} z^2 + \frac{1}{6}z^3 + \frac{1}{24}z^4,
\end{align}
where $z = -u \Delta t (i\hat{k} - \ep \hat{k}^q) + \Delta t \gamma i^{2k} \hat{k}^{2k}$. Requiring $|G(z)| \leq 1$ allows us to compute $\gamma$ for all RK methods.  We demonstrate the calculation of $\gamma$ for RK4. First, partition $z = -u \Delta t (i\hat{k} - \ep \hat{k}^q) + \Delta t \gamma i^{2k} \hat{k}^{2k}$ as $z = z_1 + z_2$, where
\begin{align}
z_1 = -u \Delta t i \hat{k}, \ z_2 = u \Delta t \ep \hat{k}^q + \Delta t \gamma i^{2k} \hat{k}^{2k}.
\end{align}
In this form, it is clear that $z_1$ arises naturally from the advection equation, while $z_2$ arises from both the spurious growth mode and the hyperviscosity term. Substituting $z = z_1 + z_2$ into the expression for $G(z)$ for RK4, we obtain
\begin{align}
G(z) = 1 + z_1 + z_2 + \frac{1}{2}\lf(z_1 + z_2\rt)^2 + \frac{1}{6}\lf(z_1 + z_2\rt)^3 + \frac{1}{24} \lf(z_1 + z_2\rt)^4.
\end{align}
Now, partition $G(z)$ as $G(z) = G_1 (z_1) + G_2(z_1,z_2)$, where
\begin{align}
G_1(z_1) &= 1 + z_1 + \frac{1}{2} z_1^2 + \frac{1}{6}z_1^3 + \frac{1}{24}z_1^4, \\
G_2(z_1,z_2) &= z_2\lf[\frac{{z_1}^3}{6} + \frac{{z_1}^2}{2} + z_1 + 1 + z_2\left(\frac{{z_1}^2}{4} + \frac{z_1}{2} + \frac{1}{2}\right) + z_2^2\left(\frac{z_1}{6} + \frac{1}{6}\right) + \frac{z_2^3}{24}\rt].
\end{align}
In the absence of growth and hyperviscosity terms, stability would require $|G_1(z)| \leq 1$. If $|G(z)|\leq 1$ and $|G_1(z_1)|\leq 1$, it follows that we require $G_2(z_1,z_2) \leq 0$. Since $G_2$ may not even be real, the easiest way to compute $\gamma$ in practice is to enforce $G_2(z_1,z_2) = 0$. In other words, we require that
\begin{align}
z_2\lf[\frac{{z_1}^3}{6} + \frac{{z_1}^2}{2} + z_1 + 1 + z_2\left(\frac{{z_1}^2}{4} + \frac{z_1}{2} + \frac{1}{2}\right) + z_2^2\left(\frac{z_1}{6} + \frac{1}{6}\right) + \frac{z_2^3}{24}\rt] &= 0. \label{eq:long_rk4}
\end{align}
To satisfy this, we set $z_2 = 0$, \emph{i.e.},
\begin{align}
u \Delta t \ep \hat{k}^q + \Delta t \gamma i^{2k} \hat{k}^{2k} = 0.
\label{eq:temp}
\end{align}
Solving \eqref{eq:temp} for $\gamma$, we get
\begin{align}
\gamma = (-1)^{1-k} \hat{k}^{q-2k}u\ep,
\end{align}
The same expression for $\gamma$ is obtained in the case of all explicit RK methods (not shown). The factor $\ep$, which is the spurious growth mode, can be estimated numerically. We also know that the maximum wavenumber $\hat{k}$ that can be resolved on a regular grid of cell width $h$ is $2h^{-1}$ (for any spatial dimension $d$). Substituting $2h^{-1}$ in place of $\hat{k}$ gives us
\begin{align}
\gamma = (-1)^{1-k}  2^{q-2k} h^{2k-q}u \ep.
\label{eq:gamma_exp}
\end{align}
Note that the above expression for $\gamma$ has an explicit dependence on the velocity $u$, the largest real part $\ep$, and the growth exponent $q$. Further, the formula contains a factor of $2^{q-2k}$, much like the $2^{-6}$ factor that appears in~\cite{FlyerNS,BarnettPHS}. We defer a discussion on $q$ to a later section.

\subsubsection{Time-stepping with explicit multistep methods}

We next calculate $\gamma$ for the case of explicit multistep methods. To illustrate the procedure, consider the Adams-Bashforth formula of order 2 (AB2)~\cite{Ascher97} applied to the generic ODE
\begin{align}
\frac{\partial c}{\partial t} = \mu c.
\end{align}
The AB2 discretization of this ODE is
\begin{align}
c^{n+1} = c^n + \frac{3}{2}\Delta t \mu c^n - \frac{1}{2}\Delta t \mu c^{n-1}.
\end{align}
Now, let $c^{n-1} = \hat{c}\lf(t_{n-1}\rt) e^{i\hat{k}x}$ so that $c^n = G c^{n-1}$ and $c^{n+1} = G^2 c^{n-1}$,where $G$ is the growth factor as before.  Substituting these equations into the AB2 scheme, we obtain the following quadratic equation for the growth factor $G$:
\begin{align}
G^2 - G\lf(1 + \frac{3}{2}\Delta t \mu\rt) + \frac{1}{2} \Delta t \mu = 0.
\label{eq:gf_ab2_1}
\end{align}
If AB2 is applied to \eqref{eq:a_hyp_1D_aux}, it follows from \eqref{eq:gf_ab2_1} that
\begin{align}
G^2 - G\lf(1 + \frac{3}{2}z_1 + \frac{3}{2} z_2\rt) + \frac{1}{2}z1 + \frac{1}{2} z_2 = 0,
\end{align}
where $z_1 = -u\Delta t i \hat{k}$ and $z_2 = u \Delta t \ep \hat{k}^q + \Delta t \gamma i^{2k} \hat{k}^{2k}$. This equation can be rewritten as
\begin{align}
\underbrace{G^2 - G\lf(1 + \frac{3}{2}z_1\rt) + \frac{1}{2}z_1}_{t_1} + \underbrace{z_2\lf(\frac{1}{2} - \frac{3}{2}G \rt)}_{t_2} = 0.
\end{align}
The term $t_1$ is free of spurious growth (and hyperviscosity), and is the term obtained by applying the AB2 method to the advection equation. If $t_2=0$ and a stability criterion on the term $u \Delta t \hat{k}$ is met, we have $|G|\leq 1$. A straightforward way to obtain an expression for $\gamma$ is therefore to require $t_2$ to vanish entirely. More specifically, we require
\begin{align}
z_2\lf(\frac{1}{2} - \frac{3}{2}G\rt) = 0.
\end{align}
Since $G$ may take on a wide range of values despite the constraint $|G|\leq 1$, the only way to guarantee the above equality is to require $z_2 = 0$. This gives us
\begin{align}
u\Delta t \ep \hat{k}^q + \Delta t \gamma i^{2k} \hat{k}^{2k} = 0,
\end{align}
which yields
\begin{align}
\gamma &= (-1)^{1-k} 2^{q-2k} h^{2k-q} u \ep,
\end{align}
which is identical to the formula for $\gamma$ in the RK4 case. This analysis is easily done for all explicit multistep methods, and even for the explicit parts of IMEX multistep methods. For an example of the latter, consider the explicit/extrapolated BDF4 method (the explicit part of the IMEX-BDF4/SBDF4 method~\cite{Ascher97}) applied to the same ODE:
\begin{align}
\frac{\frac{25}{12}c^{n+1} - 4c^n + 3c^{n-1} - \frac{4}{3}c^{n-2} + \frac{1}{4}c^{n-3}}{\Delta t} = \mu\lf(4 c^n - 6c^{n-1} + 4c^{n-2} - c^{n-3}\rt).
\end{align}
The corresponding equation for $G(z)$ is now quartic:
\begin{align}
G^4 - \frac{48}{25}G^3(1 + z) + \frac{36}{25} G^2 (1 + 2z) - \frac{16}{25}G(1 + 3z) + \frac{3}{25}(1 + 4z) = 0.
\end{align}
Writing $z = z_1 + z_2$ as before and splitting the equation into two terms $t_1(z_1)$ and $t_2(z_2,G)$ yields
\begin{align}
t_1 &= G^4 - \frac{48}{25}G^3(1 + z_1) + \frac{36}{25} G^2 (1 + 2z_1) - \frac{16}{25}G(1 + 3z_1) + \frac{3}{25}(1 + 4z_1), \\
t_2 &= z_2\lf(-\frac{48}{25}G^3 + \frac{72}{25}G^2 - \frac{48}{25}G + \frac{12}{25} \rt), \\
t_1 &+ t_2 = 0.
\end{align}
To have $t_2 = 0$ for stability, we again require $z_2 = 0$. This yields the same expression for $\gamma$ as in the AB2 and RK4 cases. As an aside, note that it is possible to obtain either $G_2 = 0$ in the RK case or $t_2 = 0$ in the multistep case if the term $z_1$ or the growth factor $G$ satisfies certain constraints (in addition to the CFL constraint). However, these cases are very rare and obtained only when $\Delta t$, $h$, and $\ep$ take on specific values.

\subsection{Selecting $\gamma$ in the diffusion-dominated regimed ($u \to 0$)}
\label{sec:gamma_2}

Another important limit is the diffusive regime where $\nu$ dominates over $u$. Current RBF-FD methods are not guaranteed to produce purely negative real eigenvalues when approximating the Laplacian $\Delta$. While a mild spread along the imaginary axis is acceptable for parabolic problems, eigenvalues with positive real parts are undesirable as they would cause spurious solution growth. To calculate the amount of hyperviscosity-based stabilization required to counter this growth, we use the same approach as in Section \ref{sec:gamma_1}, and define the auxiliary 1D PDE
\begin{align}
\frac{\partial c}{\partial t} = \nu \frac{\tilde{\partial}^2 c}{\partial x^2} + \gamma \frac{\partial^{2k}c}{\partial x^{2k}}, k \in \mathbb{N},
\label{eq:diff_hyp}
\end{align}
where the operator $\frac{\tilde{\partial}^2} {\partial x^2}$ acting on a plane wave produces the correct decay mode as well as a spurious growth mode, so that if $c(x,t) = \hat{c}(t) e^{i\hat{k}x}$, $\frac{\tilde{\partial}^2 c}{\partial x^2} = \lf(-\hat{k}^2 + \eta \hat{k}^q\rt) \hat{c}(t) e^{i\hat{k}x}$, where $\eta \geq 0$ is responsible for the spurious growth. We now discretize \eqref{eq:diff_hyp} in time using the second-order backward difference formula (BDF2) (which is A-stable)~\cite{Ascher97}, obtaining
\begin{align}
\frac{3c^{n+1} - 4c^n + c^{n-1}}{2\Delta t} = \nu \frac{\tilde{\partial}^2 c^{n+1}} {\partial x^2} + \gamma \frac{\partial^{2k}}{\partial x^{2k}} c^{n+1}, 
\end{align}
for which the growth factor $G$ satisfies
\begin{align}
\lf(1 - \frac{2}{3}\Delta t \mu\rt)G^2 -\frac{4}{3}G + \frac{1}{3} = 0,
\label{eq:visc_hyp_1}
\end{align}
where $\mu = \nu\lf(-\hat{k}^2 + \eta \hat{k}^q\rt) + \gamma i^{2k} \hat{k}^{2k}$. Collecting all diffusion terms into term $t_1$, and spurious growth and hyperviscosity terms in term $t_2$gives:
\begin{align}
&t_1 = \lf(1 + \frac{2}{3}\Delta t \hat{k}^2\rt)G^2 - \frac{4}{3}G + \frac{1}{3}, \ t_2 = -\frac{2}{3}\Delta t \lf(\nu\eta \hat{k}^q + \gamma i^{2k} \hat{k}^{2k}\rt), \\
&t_1 + t_2 =0.
\end{align}
For stability, we require that $t_2 = 0$. This gives:
\begin{align}
\nu\eta \hat{k}^q &+ \gamma i^{2k} \hat{k}^{2k} = 0, \\
\implies \gamma &= (-1)^{1-k} 2^{q-2k} h^{2k-q} \nu \eta.
\label{eq:gamma_visc_1}
\end{align}
If we had lumped the diffusion term into term $t_2$, we would have obtained an expression for $\gamma$ that accounts for the stabilization afforded by natural diffusion. We do not take this approach in this article.

\subsection{Selecting $\gamma$ for coexisting advection and diffusion ($u\neq 0$, $\nu \neq 0$)}
\label{sec:gamma_3}

Having examined both the advective and diffusive limits for reasonable choices of time integrators, we now consider the important case where both advection and diffusion may contribute significantly to the transport process. In this scenario, it is common to treat the diffusion term implicitly in time, and the advection term explicitly. Consequently, it is reasonable (in terms of computational efficiency) to also treat the hyperviscosity term implicitly in time. The relevant auxiliary PDE with hyperviscosity is now
\begin{align}
\frac{\partial c}{\partial t} + u \frac{\tilde{\partial} c}{\partial x}=  \nu \frac{\tilde{\partial}^2 c}{\partial x^2} + \gamma \frac{\partial^{2k}}{\partial x^{2k}}.
\label{eq:ad_1d}
\end{align}
For plane waves, this yields the ODE
\begin{align}
\frac{\partial c}{\partial t} = \mu_1 c + \mu_2 c,
\label{eq:ode2}
\end{align}
where
\begin{align}
\mu_1 = -u i\hat{k} + u\ep\hat{k}^{q_1}, \ \mu_2 = \nu\lf(-\hat{k}^2 + \eta \hat{k}^{q_2}\rt) + \gamma i^{2k} \hat{k}^{2k},
\end{align}
where $q_1$ and $q_2$ are the growth exponents for the approximate gradient and laplacian respectively. To keep our analysis simple, we discretize \eqref{eq:ode2} using the IMEX-BDF2 or SBDF2 method~\cite{Ascher97}; however, our technique carries over to higher order SBDF methods. The SBDF2 discretization is given by
\begin{align}
\frac{3c^{n+1} - 4c^n + c^{n-1}}{2\Delta t} = 2\mu_1 c^n - \mu_1 c^{n-1} + \mu_2 c^{n+1}.
\end{align}
Substituting in a plane wave expression yields the following equation for the growth factor $G$:
\begin{align}
\lf(1 - \frac{2}{3}\mu_2 \Delta t\rt) G^2 - \frac{4}{3}G\lf(1 + \Delta t \mu_1\rt) + \frac{1}{3}\lf(1 + 2 \Delta t \mu_1\rt) = 0.
\end{align}
Collecting advection and diffusion terms into the term $t_1$, and spurious growth and hyperviscosity terms into the term $t_2$, we obtain
\begin{align}
t_1 &= \lf(1 + \frac{2}{3}\nu\hat{k}^2 \Delta t\rt) G^2 - \frac{4}{3}G\lf(1 - \Delta t u i \hat{k} \rt) + \frac{1}{3}\lf(1 - 2 \Delta t u i \hat{k}\rt), \\
t_2 &= -\frac{2}{3}\Delta t G^2 \lf(\nu\eta\hat{k}^{q_2} + \gamma i^{2k} \hat{k}^{2k}\rt) + \frac{2}{3}\Delta t u \ep \hat{k}^{q_1} \lf(1 - 2G\rt), \\
t_1 &+ t_2 = 0.
\end{align}
As before, we require that $t_2 = 0$, \emph{i.e.}, 
\begin{align}
-G^2 \lf(\nu\eta\hat{k}^{q_2} + \gamma i^{2k} \hat{k}^{2k}\rt) + u \ep \hat{k}^{q_1} \lf(1 - 2G\rt) = 0.
\end{align}
Now, let $\gamma = \gamma_1 + \gamma_2$. Then, we have
\begin{align}
-G^2 \lf(\nu\eta\hat{k}^{q_2} + \gamma_1 i^{2k} \hat{k}^{2k} + \gamma_2 i^{2k} \hat{k}^{2k}\rt) + u \ep \hat{k}^{q_1} \lf(1 - 2G\rt) = 0.
\end{align}
We can select $\gamma_2$ to cancel out the growth term $\eta \hat{k}^{q_2}$, giving us
\begin{align}
\gamma_2 = (-1)^{1-k} 2^{q_2 - 2k} h^{2k - q_2} \eta \nu,
\end{align}
and leaving us with
\begin{align}
-G^2 \gamma_1 i^{2k} \hat{k}^{2k} + u \ep \hat{k}^{q_1} \lf(1 - 2G\rt) = 0.
\end{align}
Dividing the above equation by the non-zero quantity $-\gamma_1 i^{2k} \hat{k}^{2k}$, we get the quadratic equation
\begin{align}
G^2  + \frac{u \ep \hat{k}^{q_1}}{\gamma_1 i^{2k} \hat{k}^{2k}} \lf(2G-1\rt) = 0.
\end{align}
The roots of this equation are given by
\begin{align}
G = -\frac{u \ep \hat{k}^{q_1}}{\gamma_1 i^{2k} \hat{k}^{2k}} \pm \frac{1}{\gamma_1 i^{2k} \hat{k}^{2k}}\sqrt{\lf(u \ep \hat{k}^{q_1}\rt)\lf(u \ep \hat{k}^{q_1} + \gamma_1 i^{2k} \hat{k}^{2k}\rt)}.
\end{align}
Stability requires that $|G|\leq 1$. This can be achieved if
\begin{align}
u \ep \hat{k}^{q_1} + \gamma_1 i^{2k} \hat{k}^{2k} = 0,
\end{align}
which gives
\begin{align}
\gamma_1 = (-1)^{1-k} 2^{q_1 - 2k} h^{2k - q_1} u \ep.
\end{align}
Finally, since $\gamma = \gamma_1 + \gamma_2$, we have
\begin{align}
\gamma &= (-1)^{1-k} 2^{q_1 - 2k} h^{2k - q_1} u \ep + (-1)^{1-k} 2^{q_2 - 2k} h^{2k - q_2} \eta \nu.
\end{align}
From the analysis, it is clear that the hyperviscosity for spurious growth in advection and diffusion is \emph{additive} under the assumptions of our model. This expression naturally encodes the other cases seen so far. For instance, in the limit $\nu \to 0$, we recover \eqref{eq:gamma_exp}, while in the limit $u \to 0$, we recover \eqref{eq:gamma_visc_1}. The expressions for $\gamma$, though derived using a 1D model, generalize straightforwardly to higher dimensions by replacing $u$ with $\|\vu\|_{\infty}$ (the maximum pointwise spatial velocity or an estimate thereof).

In this article, our simulation results are expressed in terms of the \emph{Peclet number}, which is given $\rm{Pe} = \frac{uL^*}{\nu}$, where $L^*$ is some characteristic length scale. To get a higher Peclet number, one can either fix $\nu$ and make $u$ larger, or fix $u$ and make $\nu$ smaller. The former case is interesting as it magnifies the effect of the spurious eigenvalue $\ep$, while the latter serves to reduce any stabilizing diffusion. In Section \ref{sec:results}, we focus on the case of increasing $u$ for fixed $\nu$, as this is a more effective test of our hyperviscosity formulation. However, we have observed that our hyperviscosity formulation is extremely robust and stable in the scenario of decreasing $\nu$ as well; this can be seen in Section \ref{sec:application}. In both these sections, we demonstrate that it is perfectly safe to use \eqref{eq:gamma_exp} for $\gamma$ in the context of IMEX multistep methods.

\subsection{Estimating growth exponents ($q_1$ and $q_2$)}
\label{sec:models}

Thus far, our growth models have only included a \emph{generic} spurious growth term of the form $\ep \hat{k}^q_1$ or $\eta \hat{k}^q_2$. To the best of our knowledge, the literature lacks any theoretical estimates for $q_1$ or $q_2$ on any given node set. Our goal in this section is to give a simple and cost-effective \emph{numerical} procedure for computing these values on a given node set. Without loss of generality, we proceed with a description of this procedure on a 1D node set. Let $X = \{ x_j\}_{j=1}^N$ be a set of nodes on which we obtain the discrete RBF-FD differentiation matrix $G^x$ that approximates $\frac{\partial}{\partial x}$ in the advection term. Our goal now is to estimate $q_1$ for the matrix $G^x$; the procedure for $q_2$ is identical. Consider the function $f(x) = e^{i\hat{k}x}$. Its derivative is given exactly by $g(x) = i \hat{k} e^{i\hat{k}x}$. Then, define
\begin{align}
\uf = \lf.f(x)\rt|_{X},\ \ \ug = \lf.g(x)\rt|_{X} = i\hat{k}\uf,
\end{align}
where $\uf,\ug$ are vectors with $N$ entries. On the other hand, the \emph{approximate} derivative of $f(x)$ on the node set $X$ is given by
\begin{align}
\tilde{\ug} = G^x \uf.
\end{align}
According to our growth model, we know that $G^x$ is represented by the auxiliary differential operator $\frac{\tilde{\partial}}{\partial x}$. Consequently, ignoring truncation errors in $G^x$, the growth model gives
\begin{align}
\tilde{\ug} = \lf(i\hat{k} - \ep \hat{k}^{q_1}\rt) \uf.
\end{align}
This allows us to write
\begin{align}
\| \ug - \tilde{\ug}\| &= \ep \hat{k}^{q_1} \|\uf\|.
\end{align}
In the scenario with spurious growth modes, we know that $\ep \|\uf\| \neq 0$. Dividing by $\ep \|\uf\|$ and reversing the sides of the equation for clarity, we get
\begin{align}
\hat{k}^{q_1}  &= \frac{\| \ug - \tilde{\ug}\|}{\ep \|\uf\|}, \\
\implies q_1 &=\frac{\ln\lf(\| \ug- \tilde{\ug}\|\rt) - \ln\lf(\ep \|\uf\|\rt)}{\ln\hat{k}}.
\end{align}
As we remarked when deriving $\gamma$, the largest wavenumber $\hat{k}$ that can be represented on a grid of cell width $h$ is $2h^{-1}$. Thus, in the above formula for $q_1$, we use $\hat{k} = 2h^{-1}$, where $h$ is taken to be $h = N^{-1/d}$. To fully extend the above definition of $q_1$ to dimension $d>1$, we need only change the definition of $f(x)$ to be a product of exponentials in each spatial variable. The above procedure is repeated for each of the gradient differentiation matrices $G^x$,$G^y$, and $G^z$, and the largest $q_1$ value of the three is used. While any function $f(x)$ could have been chosen in principle, the use of plane waves allows a direct comparison to our growth model and an easy evaluation on any node set. In addition, to compute $q_2$, we use $g(x) = -\hat{k}^2 e^{i\hat{k}x}$ as our reference derivative function and $\tilde{\ug} =  \lf(-\hat{k}^2 + \eta\hat{k}^{q_2}\rt) \uf$, but the above formula remains the same (with $\eta$ in place of $\ep$). In practice though, since $\eta \approx 0$ and we use implicit time-stepping, it is typically unnecessary to compute $q_2$. 

To complete the above discussion, we note that $\ep$ and $\eta$ must be estimated to find $q_1$ and $q_2$. These factors correspond naturally to the real part of the eigenvalue with the largest real part. If a loose tolerance is used in the eigenvalue computation (as is done in this article), the estimates for $q_1$ and $q_2$ are likely to be an approximation as well. Note also that our definition allows $q_1$ and $q_2$ to be negative, though in practice we have observed that $q_1$ and $q_2$ are typically between $0$ and $1$. The approach detailed in this section \emph{eliminates} all free parameters from the expression for $\gamma$ for the cost of $d$ sparse matrix-vector multiplications for estimating $q$, and a few sparse matrix-vector multiplications for estimating $\ep$. In contrast, expressions for $\gamma$ in the literature contain factors that must be tuned to achieve stability~\cite{FlyerNS,FlyerLehto2012}. If theoretical developments should ever give \emph{a priori} estimates for $\ep$ and $\eta$, these can be substituted directly into the formula for $\gamma$.

\subsection{Spurious modes in the hyperviscosity operator}
\label{sec:hyp_growth}

Our analysis, for convenience, has neglected the fact that an RBF-FD approximation to the operator $\Delta^k$ may itself contain spurious eigenvalues that cause undesirable growth or decay (so that they always produce spurious growth when we consider $\gamma \Delta^k$ since $\textit{sign}(\gamma) = (-1)^{1-k}$). It is reasonable to wonder how the expression for $\gamma$ would change if one factored in this additional spurious mode. Fortunately, our analysis technique from the previous sections is perfectly amenable to this case as well. Let $\frac{\tilde{\partial}^{2k}}{\partial x^{2k}}$ be the \emph{auxiliary} hyperviscosity operator, defined by its action on plane waves:
\begin{align}
\frac{\tilde{\partial}^{2k} e^{i\hat{k}x}}{\partial x^{2k}} = \lf( i^{2k} \hat{k}^{2k} \pm \tau \hat{k}^{q_3}\rt) e^{i\hat{k}x},
\end{align}
where $\tau\geq 0$, but the spurious term has the opposite sign of the eigenvalues of the operator $\frac{\partial^{2k}}{\partial x^{2k}}$. Without loss of generality, we rederive $\gamma$ for explicit RK4 for this new scenario. Using our previous notation of $G(z)$ to represent the growth factor, splitting $z = z_1 + z_2$, where $z_1$ contains the advection term, and $z_2$ contains \emph{all} spurious growth terms, we again require $z_2 = 0$ to obtain stability. This yields the condition
\begin{align}
u\Delta t \ep \hat{k}^{q_1} + \Delta t \gamma \lf(i^{2k} \hat{k}^{2k} \pm \tau \hat{k}^{q_3}\rt) = 0,
\end{align}
where the symbols have their usual meanings. Re-arranging, we get
\begin{align}
\gamma &= -\frac{u \ep \hat{k}^{q_1}}{i^{2k} \hat{k}^{2k} \pm \tau \hat{k}^{q_3}},\\
\implies \gamma &= \frac{(-1)^{1-k} u \ep 2^{q_1-2k} h^{2k - q_1}}{1 \pm (-1)^k \tau 2^{q_3 - 2k} h^{2k-q_3}},
\end{align}
for $\hat{k} = 2h^{-1}$. Since the $\tau$ term is scaled by the typically very small value of $2^{q_3 - 2k}h^{2k-q_3}$, we find for the experiments in this article that it is safe to neglect the effect of any spurious modes in the RBF-FD approximation to the $\Delta^k$ operator when computing $\gamma$. 

\subsection{Selecting $k$}
\label{sec:k}

We next turn to the selection $k$, the power of the Laplacian in the term $\gamma \Delta^k$. It has already been observed in the RBF-FD literature that $k$ must be increased as the stencil size $n$ is increased~\cite{FlyerLehto2012}. However, the formula presented in~\cite{FlyerLehto2012} was heuristic and was fine-tuned to the node sets on the sphere. In more modern RBF-FD formulas on Euclidean domains~\cite{FlyerNS}, this approach appears to have been tested and discarded as unnecessary. However, in our experiments, we found that it was indeed necessary to scale $k$ with $n$ on scattered nodes in domains with irregular boundaries. 

We adopt an approach based on the spectral methods literature. Specifically, examining Ma's work on the Chebyshev-Legendre super viscosity method~\cite{MaSSV1,MaSSV2}, we see in Remark 2 below (5.21) that the order $k$ ($s$ in their notation) of the spectral superviscosity must scale as
\begin{align}
k \leq O(\ln N),
\end{align}
where $N$ is the total number of points. Our recipe for $k$ uses a modified version of this scaling law. Bearing in mind that an FD method reproduces a spectral method as the stencil size $n \to N$, this prompts us to select
\begin{align}
k = \lf\lfloor 1.5 \ln n\rt\rfloor.
\label{eq:k}
\end{align}
This formula has the effect of increasing $k$ with $n$ (as desired). The actual constant in front of the $\ln n$ term does not appear to have a great impact on accuracy. \revone{While the logarithmic scaling may arise naturally as a consequence of the non-uniform resolution of Chebyshev-Legendre spectral methods, it has benefits when used with RBF-FD approximations also. First, the law \eqref{eq:k} is easy to compute and requires no fine tuning. Second, it serves to increase $k$ slowly with $n$. This slow increase ensures that the stencil sizes and polynomial degrees needed for approximating $\Delta^k$ also grow only slowly with $n$. Together, the formulas for $\gamma$ and $k$ give us a hyperviscosity operator that produces the correct damping behavior in the $n \to N$ limit, \emph{and} vanishes in the $h \to 0$ limit, all without any free parameters.}

\subsection{RBF approximations to the hyperviscosity operator}
\label{sec:rbf_hyp}

As in the case of the other terms in \eqref{eq:adr_hyp}, we approximate the hyperviscosity term with PHS RBFs augmented with polynomials. First, the PHS RBF used for hyperviscosity must have sufficient smoothness so as to be able to apply the operator $\Delta^k$ to it. We use the PHS RBF
\begin{align}
\phi_{hyp}(r) = r^{2k+1},
\end{align}
which is the PHS RBF of minimal smoothness required. It follows that, $\Delta^k \phi_{hyp}(r) = C r$, where $C$ is a known constant involving $k$ and the spatial dimension $d$. This follows from repeated application of the formula for the Laplacian of a radial function in $d$-dimensions: $\Delta r^m = m(m - d + 2)r^{m-2}$. 

In addition to selecting the PHS RBF, we must also select the polynomial degree $\ell_{hyp}$ for approximating the hyperviscosity operator. Our approach here is very simple. For a PHS RBF of degree $2k+1$, the minimal polynomial degree required to prove unisolvency of the RBF interpolant is $k$~\cite{Fasshauer:2007}. We thus pick $\ell_{hyp} = k$. The observant reader may be concerned: after all, for even a first-order approximation of the operator $\Delta^k$, it would seem that $\ell_{hyp} = 2k$ is the correct choice. However, this ignores the fact that the operator in question is $\gamma \Delta^k$. Since $\gamma = O\lf(h^{2k-q}\rt)$ where we typically have $q = O(1)$, the scaling with $\gamma$ makes our approximation to this operator very high-order. When computing overlapped RBF-FD weights for the operator $\Delta^k$, the choice of $\ell_{hyp} = k$ also zeros out the derivatives of polynomial terms upto degree $\ell$, which in our case is \emph{all} terms.

The different parameters and basis function choices for our RBF-FD methods are summarized in Table \ref{tab:params}.
\begin{table}[h!]
\centering
\begin{tabular}{ccc} \toprule
		Parameter & Meaning & Value \\ \midrule
		$\ell$ & Polynomial degree for non-hyperviscosity terms & $\xi + \theta - 1$ \\
		$m$ & PHS degree & $\ell$ if $\ell$ is odd, $\ell - 1$ if $\ell$ is even \\
		$M$ & Number of polynomial terms & ${\ell + d \choose d}$ \\
		$n$ & Stencil size & $2M + \lf \lfloor \ln(2M)\rt\rfloor$\\
		$\delta$ & Overlap parameter &  0.7 if $\ell \leq 3$, 0.5 if $4 \leq \ell < 6$, 0.3 if $\ell \geq 6$\\
		$k$ & Order of hyperviscosity & $\lf\lfloor 1.5 \ln n\rt\rfloor$ \\
		$\gamma$ & Magnitude of hyperviscosity & $(-1)^{1-k} 2^{q_1-2k} h^{2k-q_1} \|\vu\|_{\infty} \ep + (-1)^{1-k} 2^{q_2-2k} h^{2k-q_2} \nu \eta$ \\
		$\phi_{hyp}(r)$ & PHS RBF for hyperviscosity & $r^{2k+1}$ \\
		$\ell_{hyp}(r)$ & Polynomial degree for hyperviscosity & $k$\\
		$n_{hyp}$ & Stencil size for hyperviscosity & $2 {\ell_{hyp} + d \choose d} + 1$\\ \bottomrule		
\end{tabular}
\caption{Table of parameters based on desired approximation order $\xi$, differential operator order $\theta$, node spacing $h$, velocity $\vu$, diffusion coefficient $\nu$, and dimension $d$.}
\label{tab:params}
\end{table}

\section{Ghost node formulation}
\label{sec:ghost}

\subsection{Method of lines}
\label{sec:mol}
In this section, we discuss our technique for time integration of \eqref{eq:adr_hyp}. In addition to selecting a time integrator, this requires the use of ghost nodes in the spatial discretization of the domain $\Omega$. \revone{As mentioned in Section \ref{sec:intro}, we observed that ghost nodes were required to enhance the global stability of our RBF-FD discretizations. This is not without precedent in the RBF-FD literature~\cite{FlyerNS,FlyerPHS,BarnettPHS}, though it is still an open question whether they are required for all problems~\cite{FlyerElliptic}}. We will now focus on solving \eqref{eq:adr_hyp} with the inclusion of a forcing term and boundary conditions:
\begin{align}
\frac{\partial c}{\partial t} &= f(c) -\vu \cdot \nabla c + \nu \Delta c + \gamma \Delta^k c, \vx \in \Omega, \\
\calB c &= g(\vx), \vx \in \partial \Omega,
\label{eq:adr_complete}
\end{align}
where $f$ is some forcing or reaction term. Our approach is to discretize the above equations in a method of lines formulation: we first discretize space, then discretize time.

We assume now that the node set $X$ contains the required ghost nodes; see~\cite{SFKSISC2017,fornbergflyerfast} for some ways to generate node sets and ghost nodes. Partition $X$ as $X = X_i \cup X_b \cup X_g$, where the subscripts $i$,$b$, and $g$ indicate interior, boundary and ghost nodes, respectively. Let $C$ be the numerical solution to \eqref{eq:adr_hyp}; write the vector $C$ as $[C_i,C_b,C_g]^T$. Let $L$ be the discrete Laplacian, $G^x$, $G^y$, and $G^z$ the discrete components of the gradient operator, and $H$ be the discrete hyperviscosity operator, all formed using overlapped RBF-FD on the node set. Let $\vu = [u_x,u_y,u_z]^T$. For each discrete differential operator $K$, we require the following matrices:
\begin{align}
K_{if} &= \begin{bmatrix} K_{ii} & K_{ib} & K_{ig} \end{bmatrix},\\
K_{bf} &= \begin{bmatrix} K_{bi} & K_{bb} & K_{bg} \end{bmatrix},
\end{align}
where $K_{ii}$ maps vectors from the interior to the interior, $K_{ib}$ from the boundary to the interior, and so forth. Finally, let $B$ be the discrete boundary condition operator approximating $\calB$. We can now write the semi-discrete analog of \eqref{eq:adr_complete} as
\begin{align}
\frac{\partial C_i}{\partial t} &= f(C_i) - (u_x)^T_i G^x_{if}C - (u_y)^T_i G^y_{if}C - (u_z)^T_i G^z_{if}C + \nu L_{if}C + \gamma H_{if} C,\label{eq:mol1} \\
\frac{\partial C_b}{\partial t} &= f(C_b) - (u_x)^T_b G^x_{bf}C - (u_y)^T_b G^y_{bf}C - (u_z)^T_b G^z_{bf}C + \nu L_{bf}C + \gamma H_{bf} C,\label{eq:mol2} \\
B_{bi} &C_i + B_{bb} C_b + B_{bb} C_g = g_b, \label{eq:mol3}
\end{align}
where $g_b = \lf.g(\vx)\rt|_{X_b}$. The above system enforces the PDE up to and including the boundary, and enforces boundary conditions at the boundary to accommodate the additional unknowns at the ghost nodes (of which there are as many as boundary nodes). Next, we proceed with the time-discretization of the constrained ODE system \eqref{eq:mol1}--\eqref{eq:mol3}. To simplify the discussion, we use forward Euler for all nonlinear terms, and backward Euler for stiff terms (including hyperviscosity). The discrete system of equations is:
{\small
\begin{align}
\frac{C_i^{n+1} - C_i^n}{\Delta t} &= f(C_i^n) - \lf((u_x)^T_i\rt)^n G^x_{if}C^n - \lf((u_y)^T_i\rt)^n G^y_{if}C^n - \lf((u_z)^T_i\rt)^n G^z_{if}C^n + \nu L_{if}C^{n+1} + \gamma H_{if} C^{n+1}, \\
\frac{C_b^{n+1} - C_b^n}{\Delta t} &= f(C_b^n) - \lf((u_x)^T_b\rt)^n G^x_{bf}C^n - \lf((u_y)^T_b\rt)^n G^y_{bf}C^n - \lf((u_z)^T_b\rt)^n G^z_{bf}C^n + \nu L_{bf}C^{n+1} + \gamma H_{bf} C^{n+1}, \\
B_{bf} &C_i^{n+1} + B_{bb} C_b^{n+1} + B_{bg} C_g^{n+1} = g_b^{n+1}.
\end{align}
}
Higher order analogues follow naturally. We use the semi-implicit backward difference formula of order 4 (SBDF4), started with a step each of SBDF1, SBDF2, and SBDF3~\cite{Ascher97}. Examining the above systems of equations, one can see that the ghost point values $C_g^{n+1}$ are a part of the system of unknowns, and that the ghost point values appear on the right hand side at time level $n$ in the first-order scheme, and time levels $n$ and $n-1$ in the second-order scheme. Clearly, except for the first time-step, there is no need to fill the ghost cells. However, the first time-step requires some careful attention in the case of a general boundary operator $\calB$.

If $\calB$ enforces Neumann or Robin conditions, the ghost values $C_g^n$ can be obtained by solving the following linear system~\cite{FlyerNS}:
\begin{align}
B_{bi} C_i^n + B_{bb} C_b^n + B_{bg} C_g^n = g_b^n, \\
\implies B_{bg} C_g^n = g_b^n - \lf(B_{bi} C_i^n + B_{bb} C_b^n\rt).
\end{align}
However, in the case of pure Dirichlet boundary conditions (as opposed to mixed, Neumann, or Robin conditions), the matrix $B_{bg}$ is a matrix of zeros and therefore has no inverse. A simple solution that works for a variety of boundary condition types is to spatially extrapolate the solution from the interior and boundary to the ghost points using local RBF interpolants of nearby nodal values. This is only done once (at time $t=0$), since all subsequent ghost node values are obtained from the time-stepping process itself.

\subsection{Positivity-preserving filter}
\label{sec:filter}

Many of our target mathematical models require that the values of the transported variable $c$ be positive at all times, with even small spurious negative values potentially leading to enormous spurious feedback loops~\cite{LEIDERMAN:2011:GWF,LEIDERMAN:2014:OMM}. High-order RBF-FD methods do exhibit a small amount of dispersion; see the figures in~\cite{FlyerNS}. While a full exploration of filters to rectify this dispersion is beyond the scope of this work, we adopt the following simple approach to ensure positivity:
\begin{enumerate}
\item Check if the elements of the solution vectors $C_i$ and $C_b$ have a negative sign.
\item Set any element with a negative sign to zero.
\end{enumerate}
This filter is applied at the end of each time-step. As we will show in the results section, this does not affect convergence rates at all in the context of advection-diffusion problems. Further, it appears to reduce errors and restore correct convergence rates on a pure advection test case. For more details, see Section \ref{sec:results}.

\section{Results}
\label{sec:results}
\begin{figure}[h!]
\centering
\subfloat[Interior and boundary nodes for the disk]
{
	\includegraphics[scale=0.4]{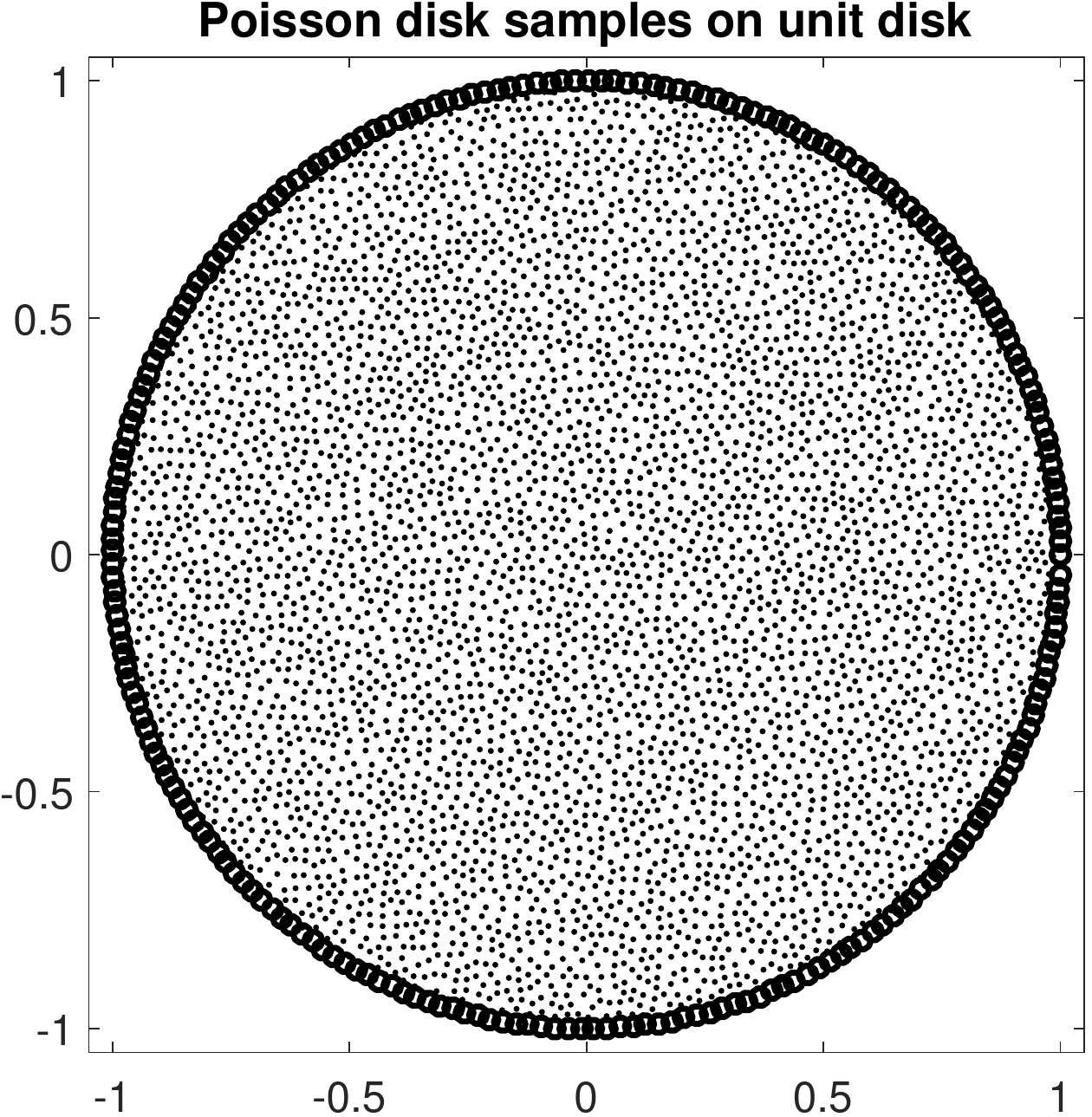}
	\label{fig:nodes2d}	
}
\subfloat[Boundary nodes for the ball]
{
	\includegraphics[scale=0.4]{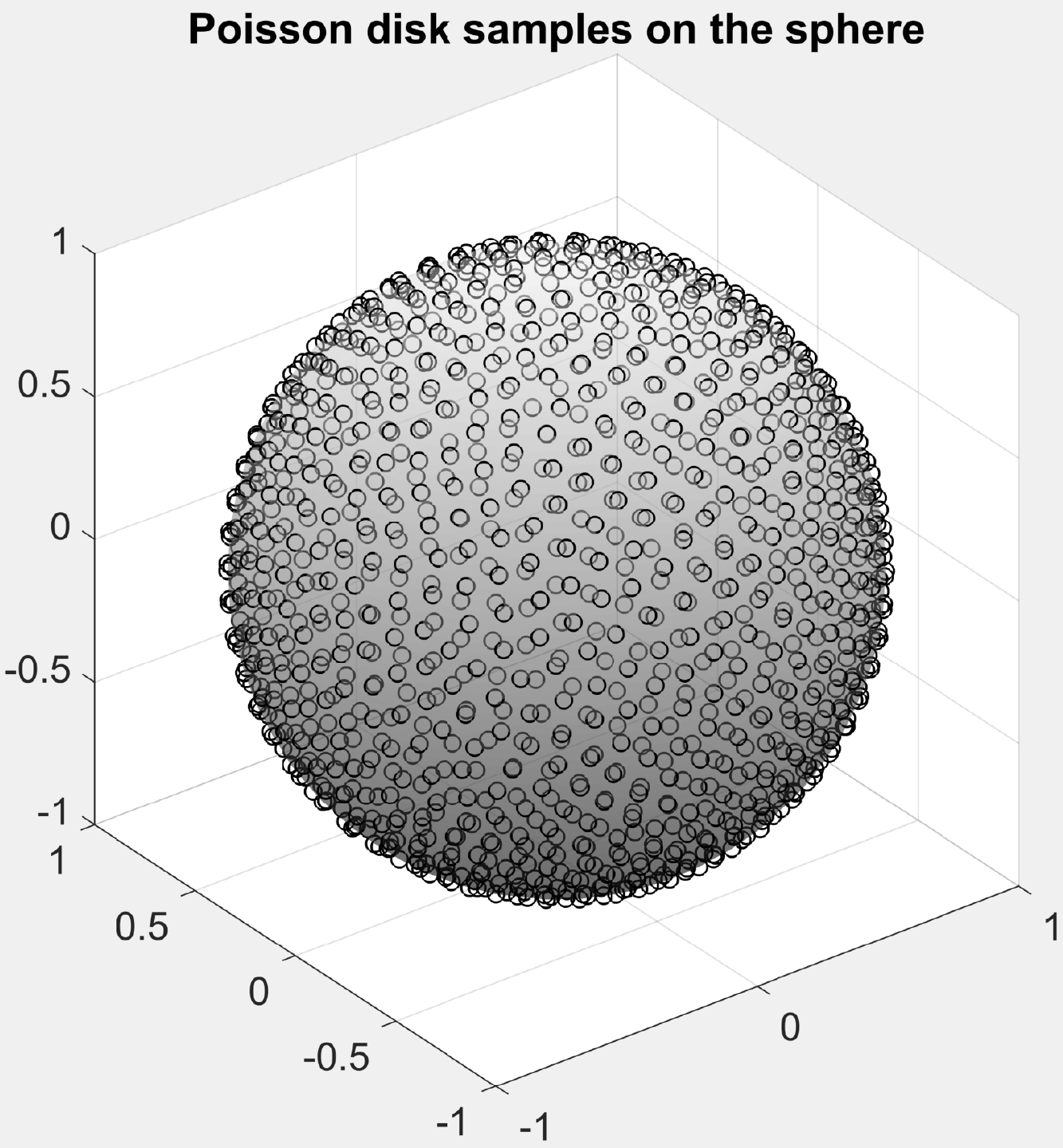}
	\label{fig:nodes3d1}	
}
\subfloat[Interior nodes for the ball]
{
	\includegraphics[scale=0.4]{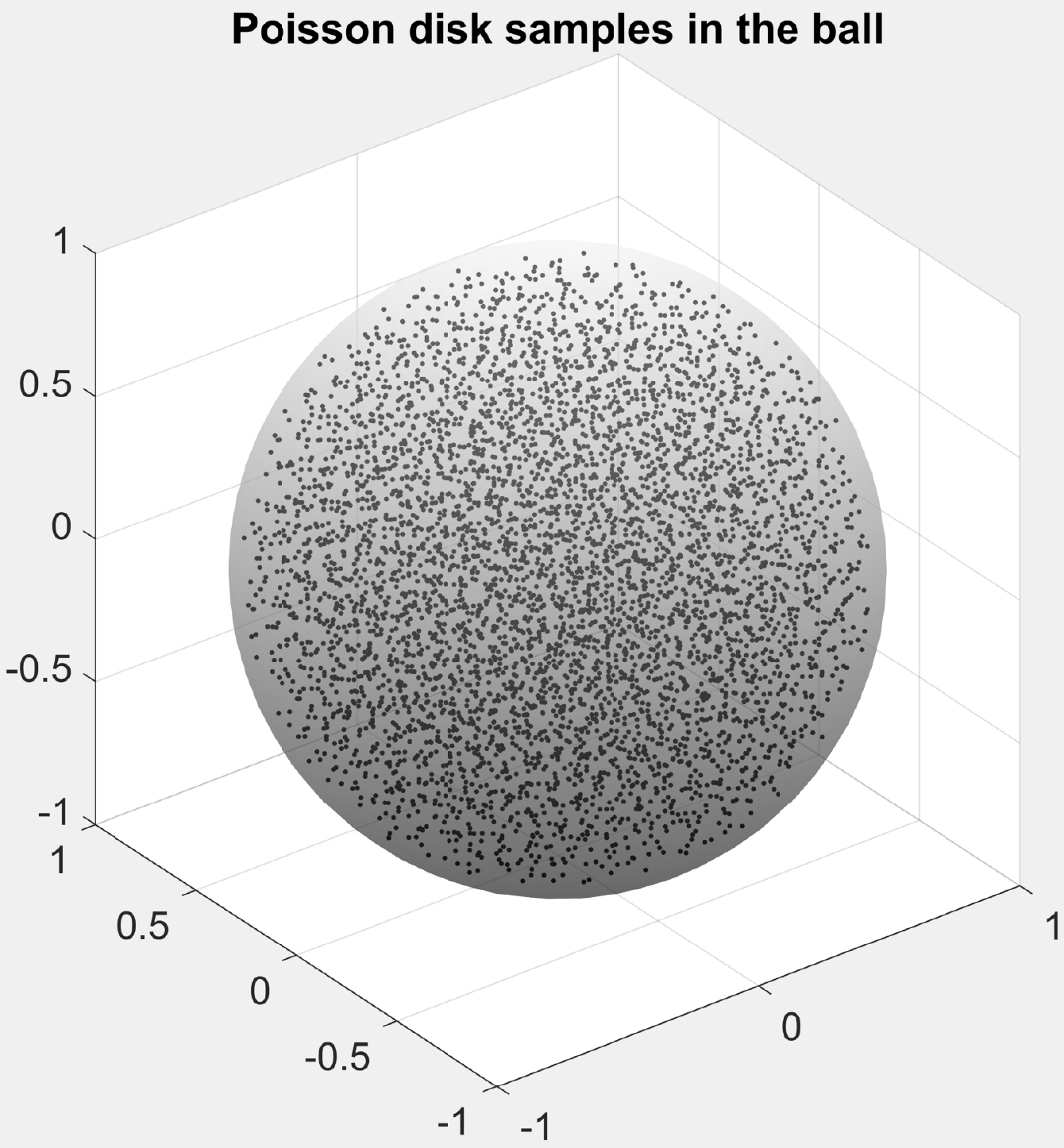}
	\label{fig:nodes3d2}	
}
\caption{Node sets on the disk and in the ball. The above node sets were generated using the algorithms in~\cite{SFKSISC2017}. Interior nodes are shown as filled circles, and boundary nodes as empty circles.}
\label{fig:nodes}
\end{figure}
We now test our hyperviscosity formulation via convergence studies on the forced advection-diffusion equation over a range of Peclet numbers. The forcing term is selected to maintain a prescribed solution for all time, and the prescribed solution is used to test spatial convergence rates. We solve this test problem on the closed unit disk in $\mathbb{R}^2$ and the closed unit ball in $\mathbb{R}^3$. Node sets were generated using the $O(N)$ node generation algorithms described in~\cite{SFKSISC2017}.  Some examples of these node sets are shown in Figure \ref{fig:nodes}. Given a true solution $c(\vx,t)$ and a numerical solution $C(\vx,t)$, we compute relative errors at the final time $t=2$ on the node set $X$ as $e_{\ell_p} = \frac{\|c_X - C_X\|_p}{\|c_X\|_p}$, where $p=2,\infty$.

\subsection{Advection on the unit disk}
\begin{figure}[h!]
\centering
\subfloat[Without filter]
{
	\includegraphics[scale=0.6]{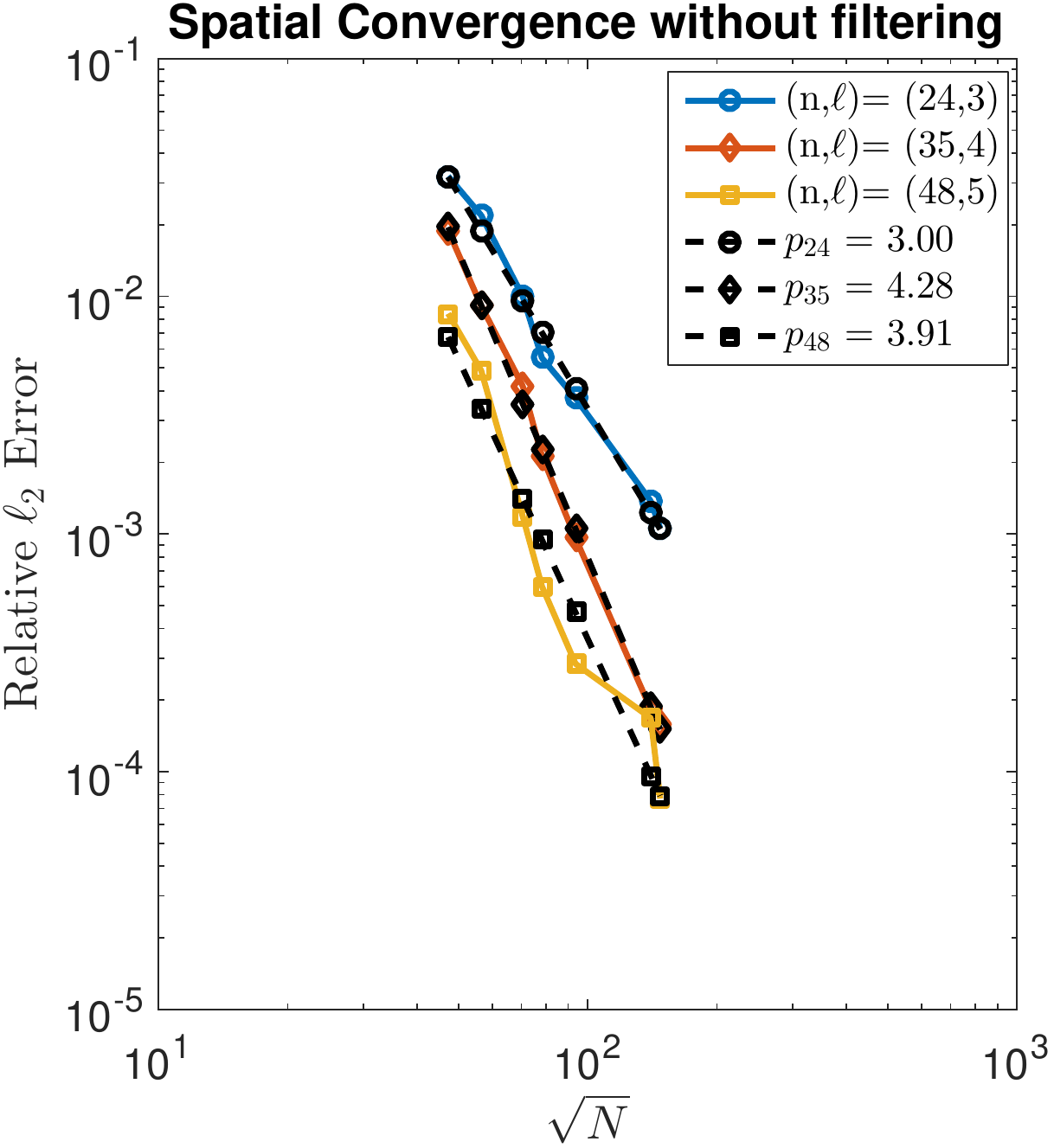}
	\label{fig:adv2d_11}	
}
\subfloat[With filter]
{
	\includegraphics[scale=0.6]{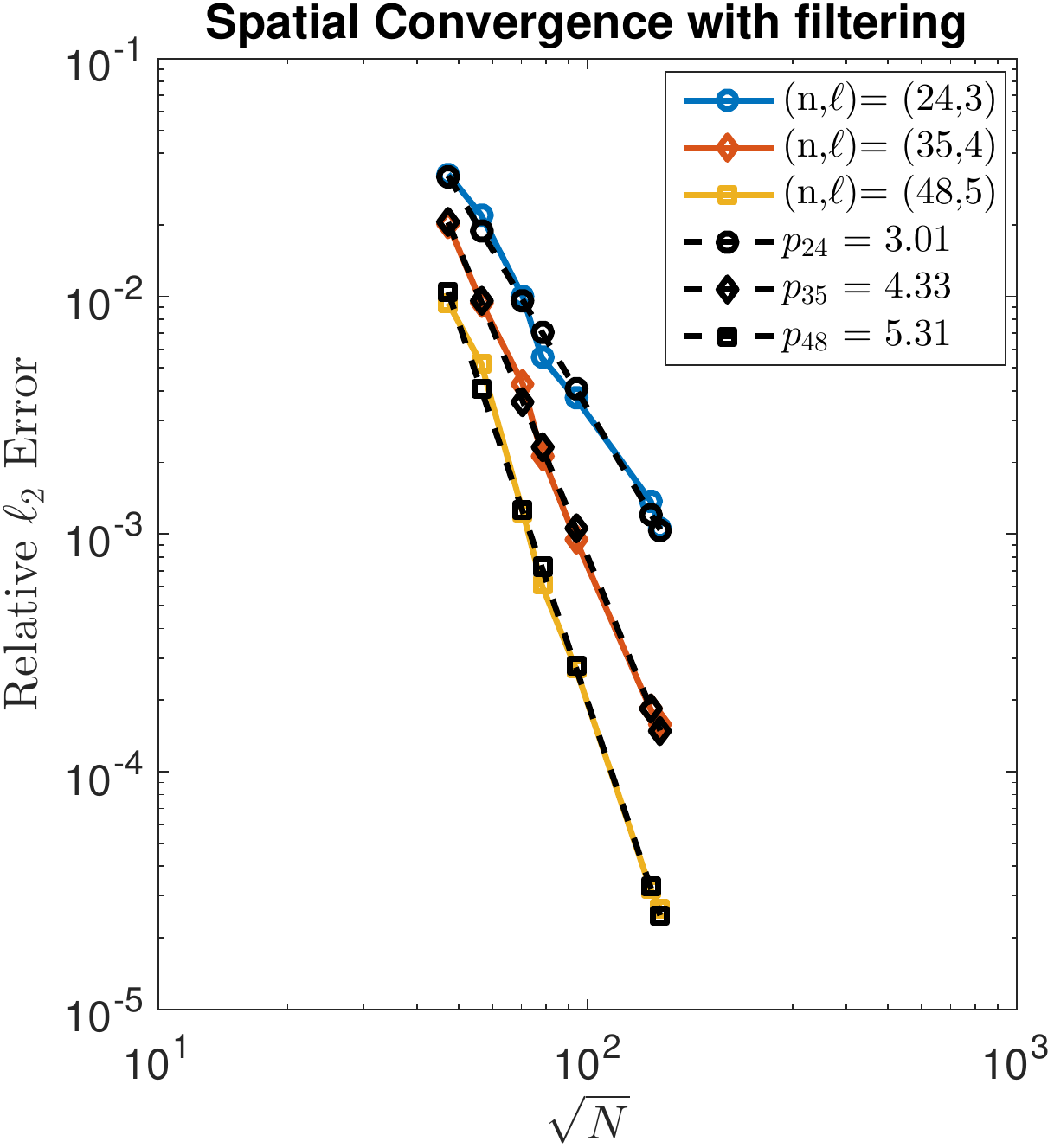}
	\label{fig:adv2d_12}	
}
\caption{Relative $\ell_2$ error vs $\sqrt{N}$ as a function of stencil size $n$ and polynomial degree $\ell$ for advection on the unit disk. The figure on the left shows results without filtering, and the figure on the right with filtering. The dashed lines are lines of best fit indicating the slope (and hence convergence rate).}
\label{fig:adv2d}
\end{figure}
Our first test involves applying our new hyperviscosity formulation to solve a pure advection problem on the unit disk. Though the focus of this article is on advection-diffusion equations, this test problem allows us to study stability and convergence rates without any stabilizing diffusion in the problem. Further, it allows us to test the hyperviscosity formulation from Section \ref{sec:hyp_exp} in the context of explicit RK methods. Finally, we also use this test to determine the effect of filtering on convergence rates.

The test involves advecting a Gaussian in a rotational velocity field. The initial condition is given by
\begin{align}
c(x,y,0) = e^{- \frac{(x - x_0)^2 + (y-y0)^2)}{\lambda^2}},
\end{align}
where $x_0 = 0.5$, $y_0 = 0.05$, and $\lambda = 1/8$. The velocity field is given by
\begin{align}
\vu(x,y,t) = 2\pi\sin(\pi t)[y,-x].
\end{align}
This velocity field effectively moves the initial condition some distance, and then brings it back to its starting position at time $t=2$, which is our chosen final time. Since the problem is a fully explicit advection equation of the form $\frac{\partial c}{\partial t} = -\vu \cdot \nabla c$, an excellent time integrator is the explicit (classical) RK4 method. In this case, since we only have a first order differential operator in the PDE ($\theta=1$), we set $\xi = \ell$. The results with and without filtering are shown in Figure \ref{fig:adv2d}. From Figure \ref{fig:adv2d_11}, we can see that increasing $\ell$ increases the convergence rate, except in the case of $\ell=5$, where we see a reduction in order of convergence. In Figure \ref{fig:adv2d_12}, we see clearly that using a filter to preserve positivity has no effect on the convergence rates when $\ell=3,4$. However, for $\ell=5$, we see that using the filter restores the correct convergence rate. Our results also demonstrate that our hyperviscosity formulation results in stable simulations in the context of RK4 timestepping, even when the hyperviscosity is stepped explicitly in time (and without any stabilizing diffusion). We use the filter in all our test results, though the errors without the filters appear to be very similar when some diffusion is present.

\subsection{Forced Advection-Diffusion on the disk}
\begin{figure}[h!]
\centering
\subfloat[]
{
	\includegraphics[scale=0.6]{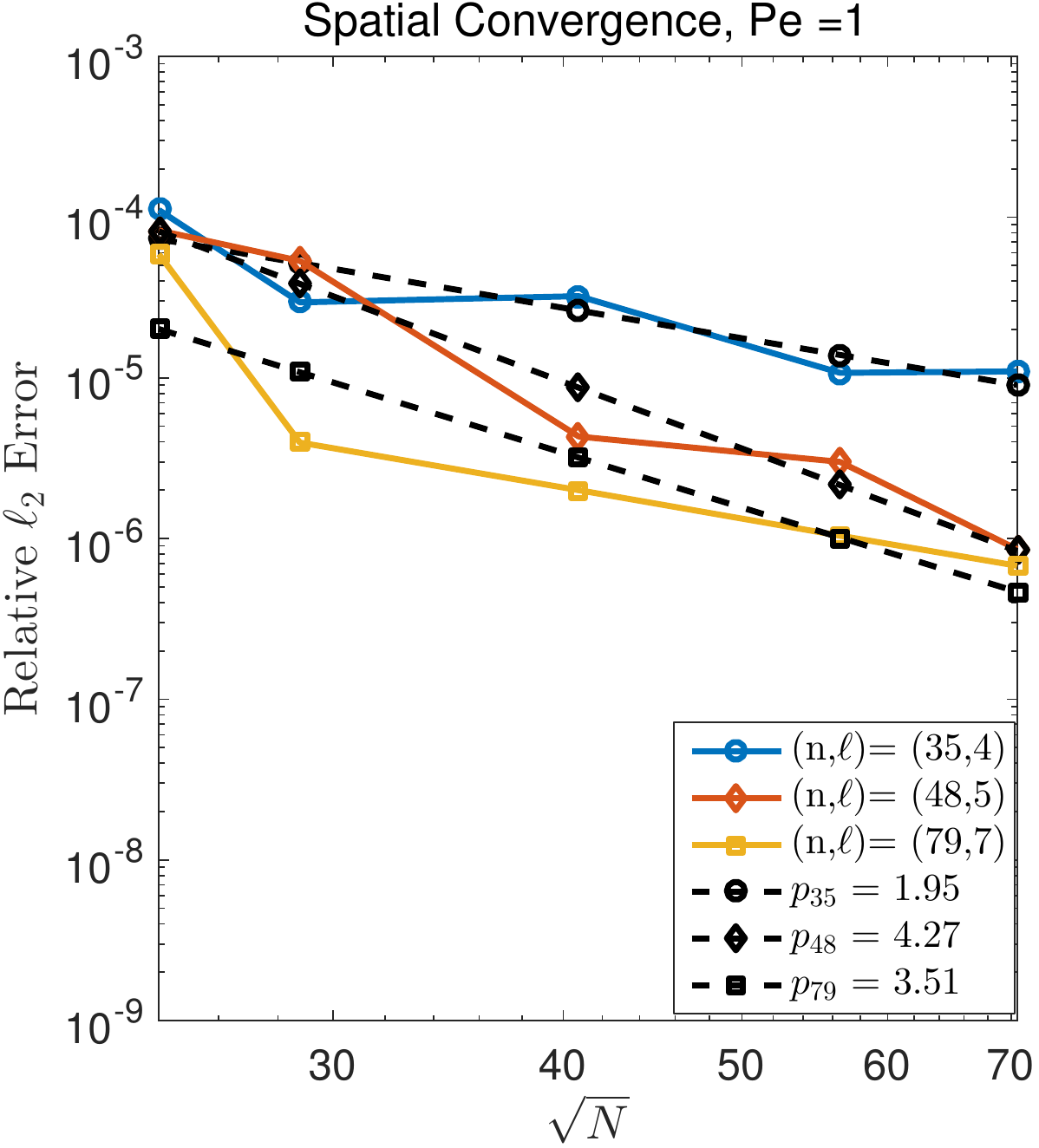}
	\label{fig:res2d_11}	
}
\subfloat[]
{
	\includegraphics[scale=0.6]{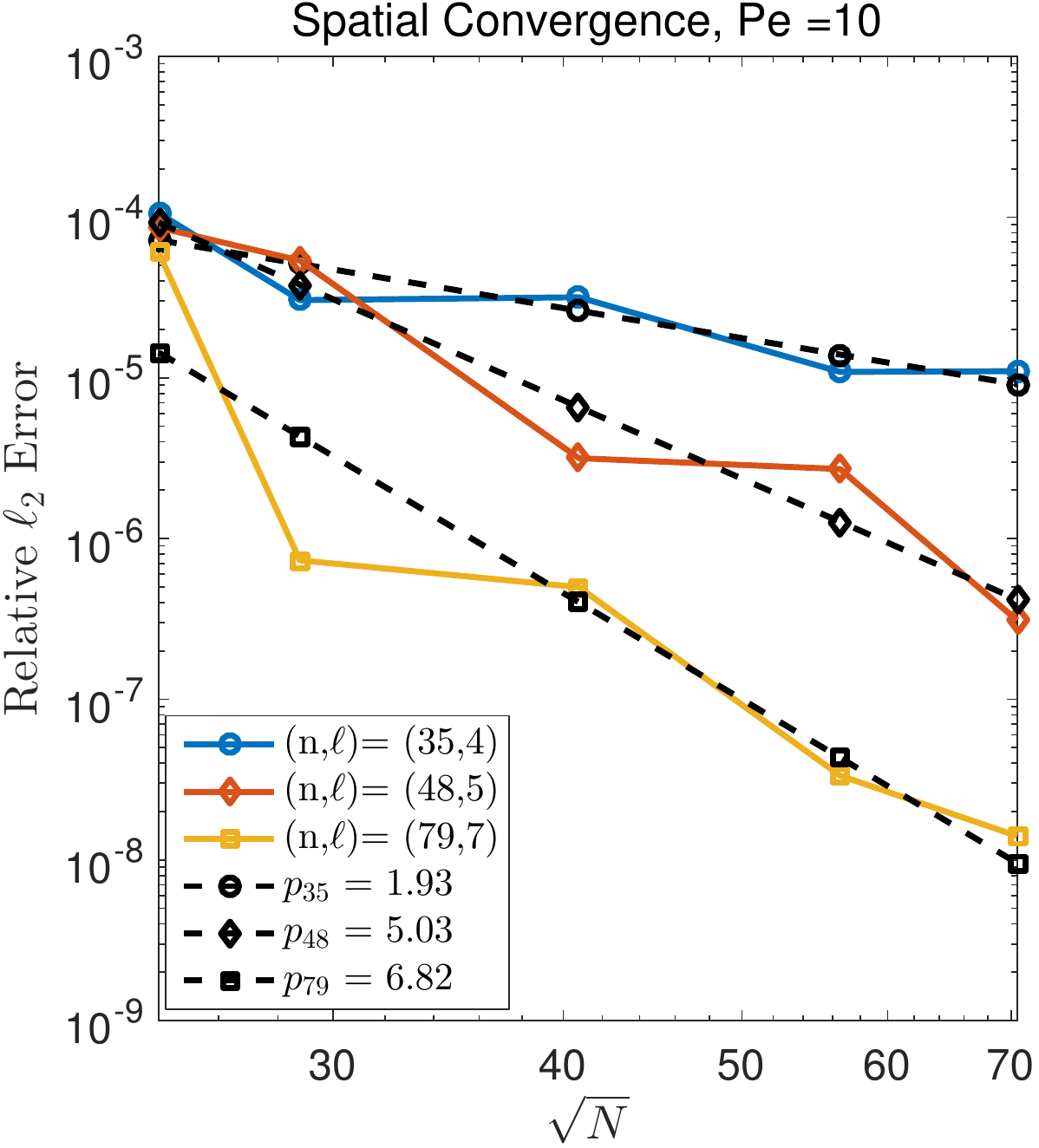}
	\label{fig:res2d_12}	
}

\subfloat[]
{
	\includegraphics[scale=0.6]{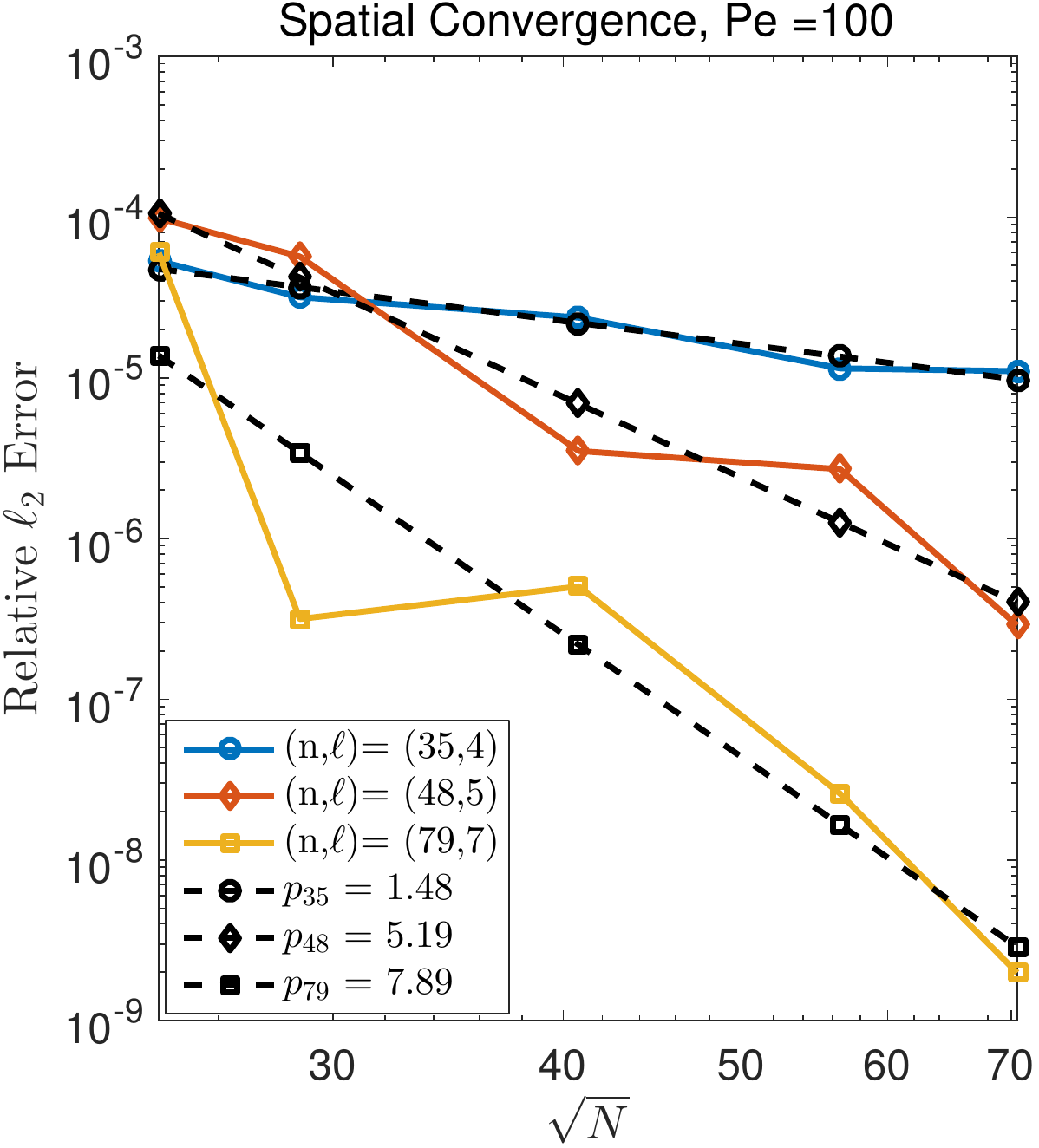}
	\label{fig:res2d_21}	
}
\subfloat[]
{
	\includegraphics[scale=0.6]{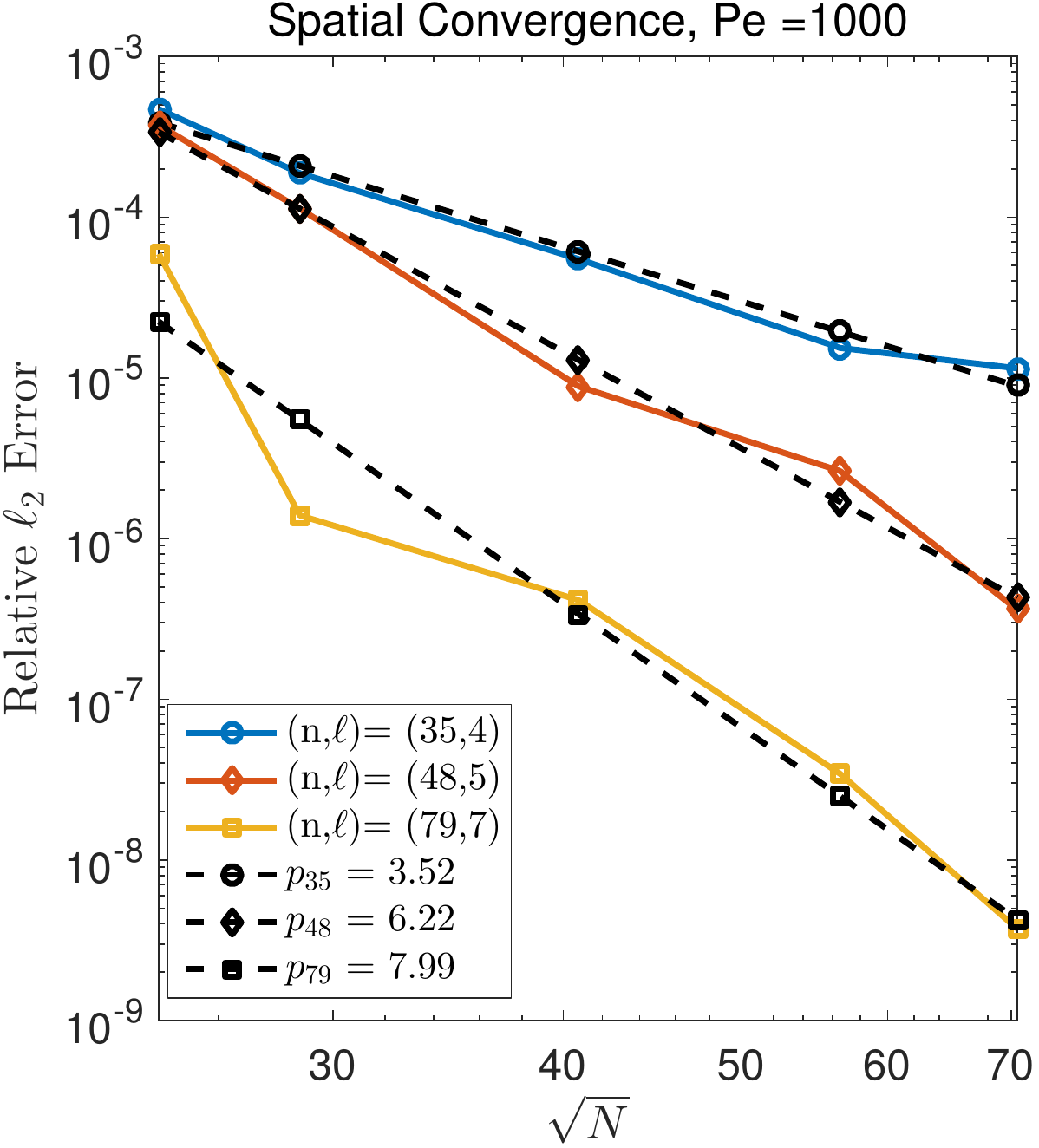}
	\label{fig:res2d_22}	
}
\caption{Relative $\ell_2$ error vs $\sqrt{N}$ as a function of stencil size $n$ and polynomial degree $\ell$ for forced advection-diffusion on the unit disk. The dashed lines are lines of best fit indicating the slope (and hence convergence rate).}
\label{fig:res2d}
\end{figure}
The next test involves solving the advection-diffusion equation on the unit disk. The manufactured solution is given by
\begin{align}
c(\vx,t) = 1 + \sin(\pi x) \cos(\pi y) e^{-\pi t},
\end{align}
and the incompressible velocity field $\vu(\vx,t) = [u,v]$ is given by
\begin{align}
\vu(\vx,t) &= \sin\lf(\pi \|\vx\|_2^2 \rt) \sin(\pi t) [y,-x].
\end{align}
This velocity field vanishes on the boundary of the disk, allowing us to safely impose Neumann boundary conditions on $c$ there using the manufactured solution. This setup allows us to measure errors in our numerical solution against the prescribed $c$. To test our hyperviscosity formulation across a range of Peclet numbers, we scale the velocity field by the Peclet number $\rm{Pe}$. In all tests shown here, we set the diffusion coefficient to $\nu = 1$. We set the time-step to 
\begin{align}
\Delta t = \frac{h}{2\|\vu\|_{max}},
\end{align}
where $\|\vu\|_{max}$ is calculated by taking the maximum of $\|\vu\|$ over space and time. We simulate the PDE to time $t=2$ using the SBDF4 method, starting up with a single step each of SBDF1, SBDF2, and SBDF3. The PDE contains both a first-order ($\theta=1$) and second-order ($\theta=2$) differential operator; rather than using different polynomial degrees $\ell$ for each operator, we use $\theta = 2$ to select $\ell$ for a given value of $\xi$ (the desired convergence rate) The results for $\ell = 4,5,7$ (with $\xi = \ell-1$) are shown in Figure \ref{fig:res2d}, plotted as a function of $\sqrt{N}$ (proportional to $1/h$).

First, it is clear from Figure \ref{fig:res2d} that increasing $\xi$ and therefore $\ell$ indeed has the effect of increasing convergence rates in general regardless of $\rm{Pe}$. Figure \ref{fig:res2d_11} shows an anomalous result for $\ell=7$ where the convergence rate drops off; however, we have verified that this is due to temporal errors dominating the total error. Indeed, at this spatial resolution and value of $\rm{Pe}$, the time-step is roughly $\Delta t = O(10^{-2})$, which incurs a fairly large error even with a fourth-order time-integrator. As $\rm{Pe}$ is increased and smaller $\Delta t$ values are consequently used, we see high order convergence rates being restored. However, as the Peclet number is increased, we see an increase in the observed order of convergence.  This is easily explained: as $\rm{Pe}$ goes up, the influence of the diffusion term (with $\theta=2$) goes down. Consequently, most of the spatial errors come from the gradient operator ($\theta=1$) since the amount of added hyperviscosity is very small. The polynomial degree $\ell$, however, was selected with $\theta=2$ in mind; it is thus reasonable to expect an extra order of convergence when $\theta=1$.
We note that similar results were seen when the magnitude of the velocity was fixed and $\nu$ was decreased to increase $\rm{Pe}$ (not shown); in those cases, as $\nu$ was decreased, it was important to scale the boundary rows of the time-stepping matrix (and the corresponding right hand sides) by $(\Delta t)^{-2}$ to ensure that the spectral radius was less than one. It appears that our proposed formula for $\gamma$ is stable across a range of parameters and Peclet numbers, as is our proposed scaling law for $m$ and $\ell$.
\subsection{Forced Advection-Diffusion in the ball}
\begin{figure}[h!]
\centering
\subfloat[]
{
	\includegraphics[scale=0.6]{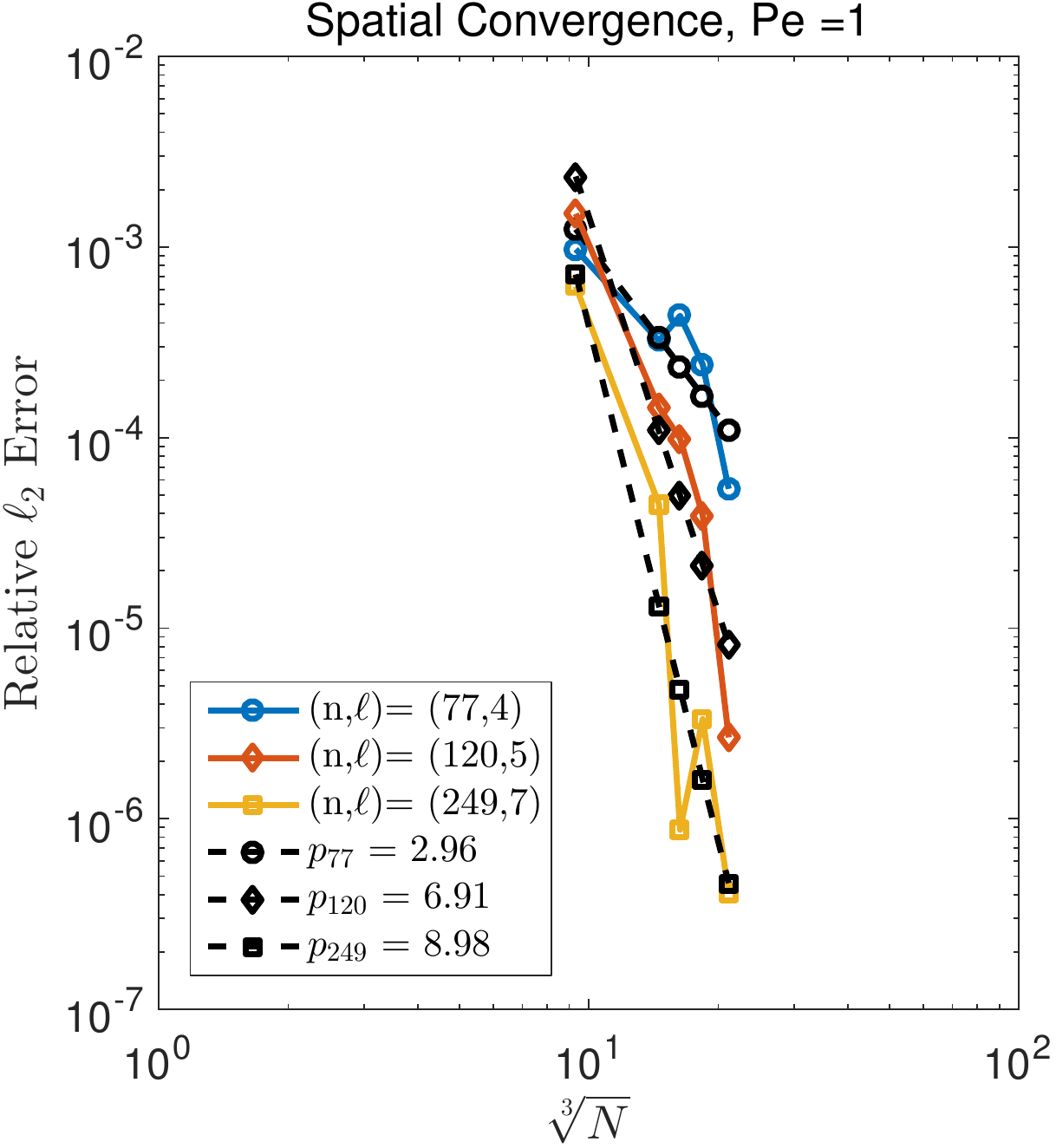}
	\label{fig:res3d_11}	
}
\subfloat[]
{
	\includegraphics[scale=0.6]{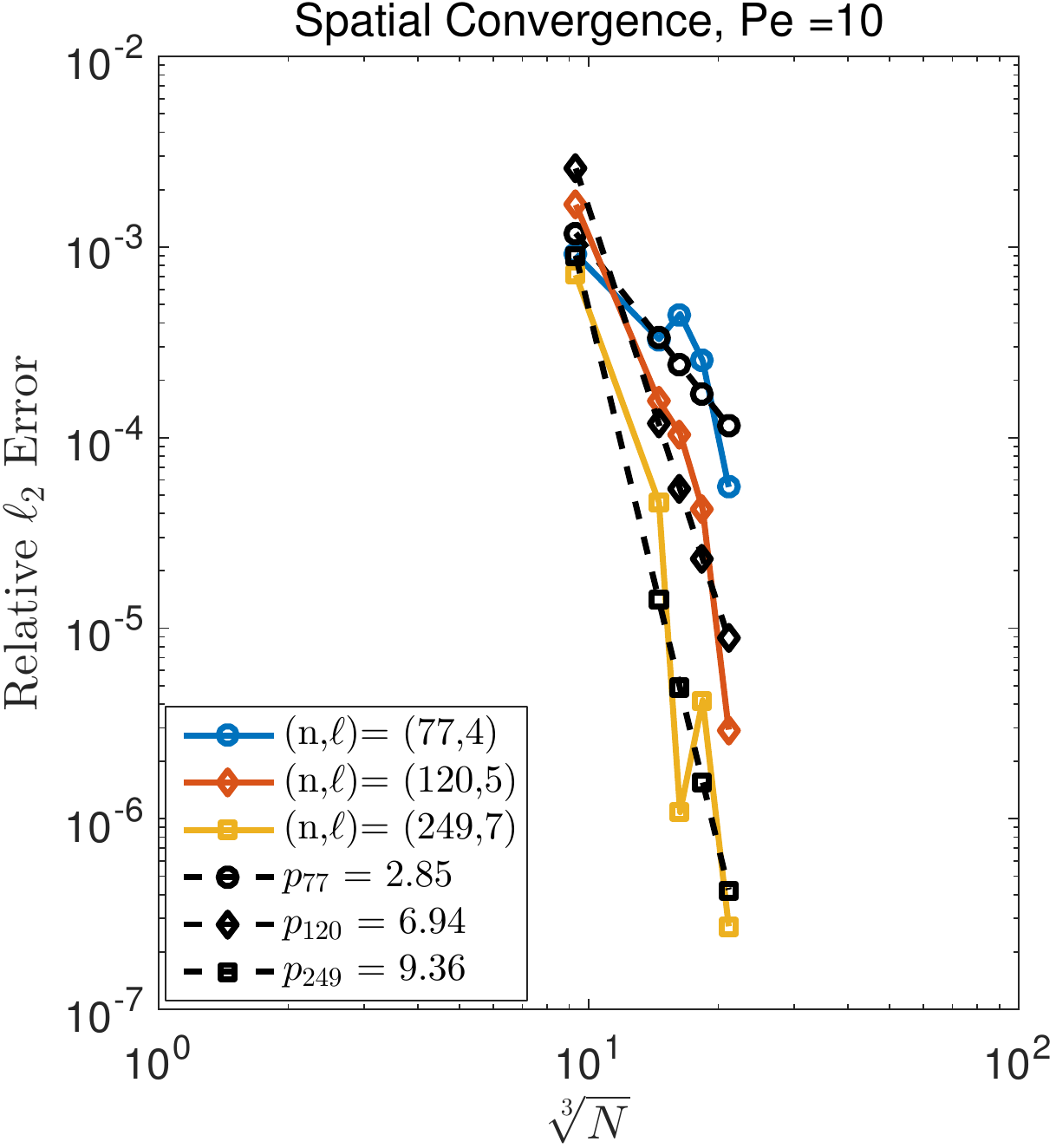}
	\label{fig:res3d_12}	
}

\subfloat[]
{
	\includegraphics[scale=0.6]{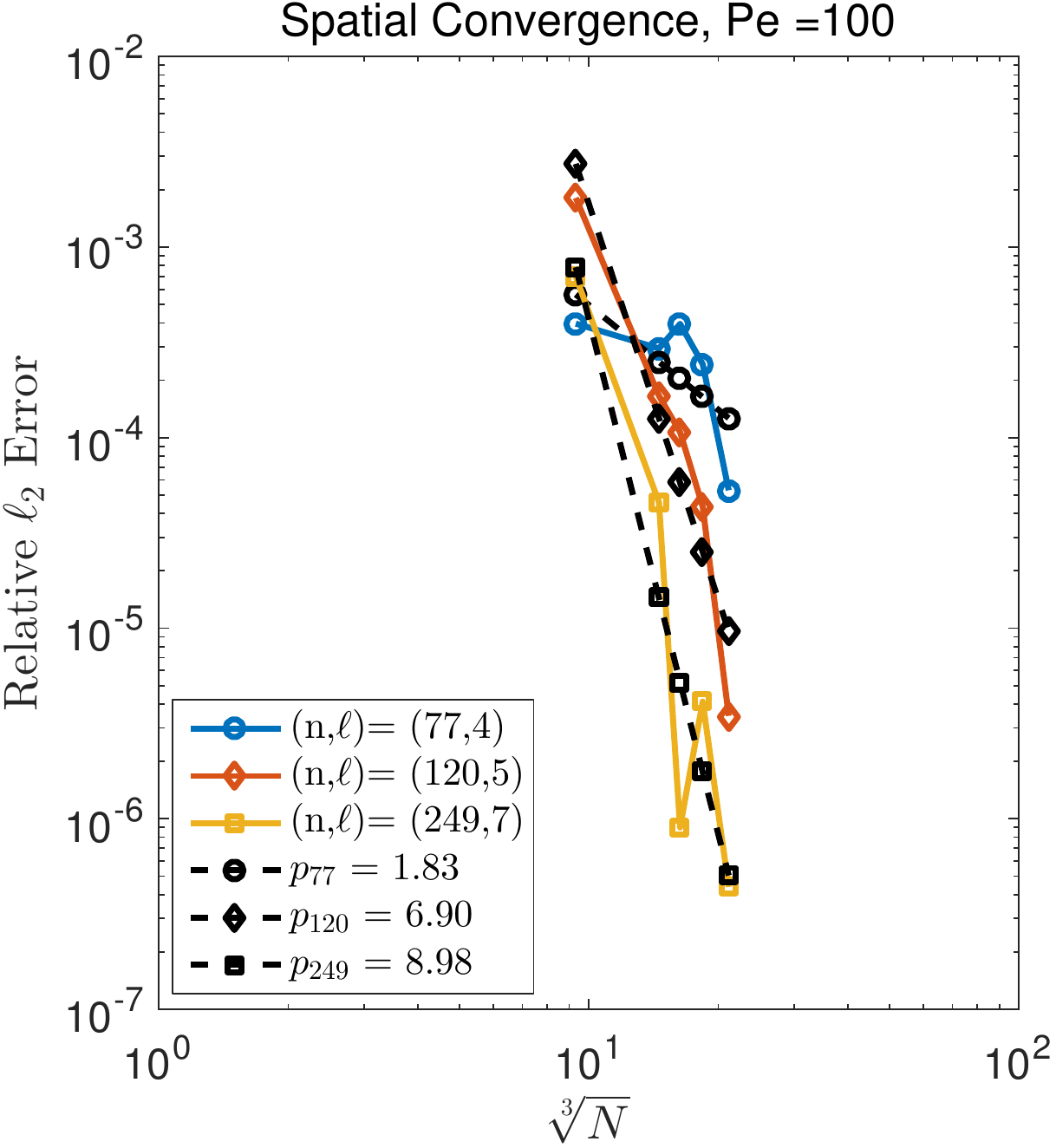}
	\label{fig:res3d_21}	
}
\subfloat[]
{
	\includegraphics[scale=0.6]{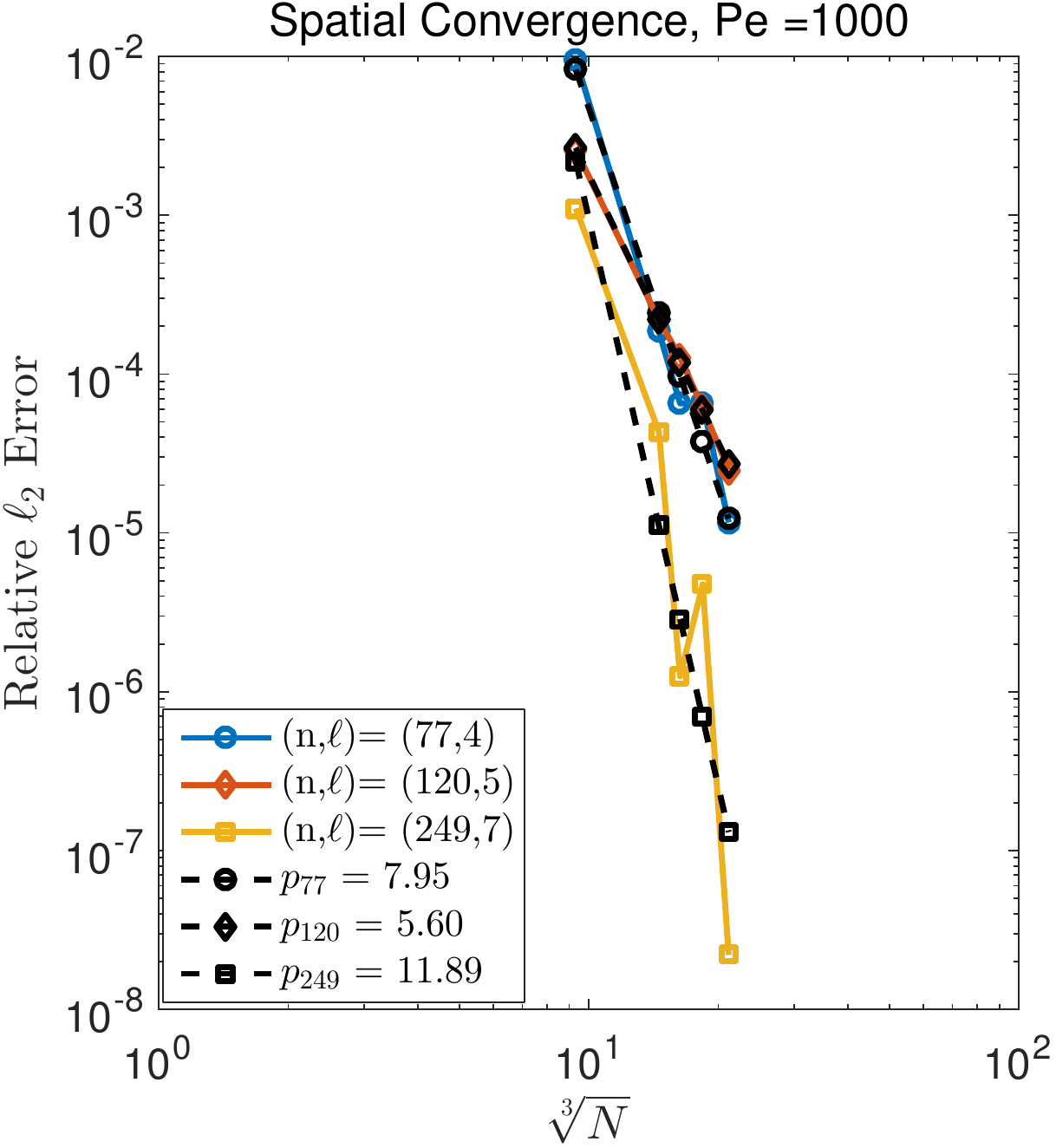}
	\label{fig:res3d_22}	
}
\caption{Relative $\ell_2$ error vs $\sqrt[3]{N}$ as a function of stencil size $n$ and polynomial degree $\ell$. The dashed lines are lines of best fit indicating the slope (and hence convergence rate).}
\label{fig:res3d}
\end{figure}
Next, we conducted the 3D analogue of the test on the disk: a convergence study on the forced advection-diffusion equation in the unit ball. In this case, our manufactured solution is
\begin{align}
c(\vx,t) = 1 + \sin(\pi x) \cos(\pi y) \sin(\pi z) e^{-\pi t},
\end{align}
and the incompressible velocity field $\vu(\vx,t) = (u,v,w)$ is given by
\begin{align}
\vu(\vx,t) = \sin\lf(\pi \|\vx\|_2^2\rt)\sin(\pi t)[-yz,2xz,-xy].
\end{align}
It is easy to show that $\vu(\vx,t)$ vanishes on the unit sphere, allowing us to impose Neumann boundary conditions on $c(\vx,t)$ there. The time-steps and Peclet numbers were chosen as in the 2D case, and we set $\nu = 1$ for all tests. We simulate the PDE to time $t=2$, and measure errors against the manufactured solution. The results are shown in Figure \ref{fig:res3d}, plotted as function of $\sqrt[3]{N}$  (proportional to $1/h$). Once again, as the Peclet number is increased, we see an increase in convergence rates, though more erratically than in the 2D case. For $\xi=3$, for instance, $\rm{Pe}=1,10$ give similar results, while $\rm{Pe}=100$ gives a reduced order of convergence and $\rm{Pe}=1000$ gives a much higher convergence rate than expected. However, in general, as $\ell$ is increased, we see an increase in convergence rates and a decrease in errors.

\section{A coupled problem in the spherical shell}
\label{sec:application}
We now apply our method to solving a coupled problem in the spherical shell (inspired by platelet aggregation and coagulation). In this problem, we track a chemical concentration $c(\vx,t)$ in a fluid inside a spherical shell domain (inner radius 0.3, outer radius 1.0). $c(\vx,t)$ is transported by the advection-diffusion equation in an incompressible velocity field, once again giving:
\begin{align}
\frac{\partial c}{\partial t} + \vu \cdot \nabla c &= \nu \Delta c + f_1(\vx,t), \vx \in \Omega,
\label{eq:coupled_c}
\end{align}
where $f_1(\vx,t)$ is some forcing term, and $\Omega$ is the region between the two spheres of radius 0.3 and 1.0. On the outer spherical boundary $\mathbb{S}_1$, $c$ satisfies a time-varying inhomogeneous Neumann boundary condition:
\begin{align}
-\nu \frac{\partial c(\vx,t)}{\partial \vn} = g(\vx,t), \vx \in \mathbb{S}_1,
\end{align}
which corresponds to an \emph{inward} flux of $c$. The inner spherical boundary $\mathbb{S}_2$ is viewed as a reactive ``zone'' on which chemicals can bind, unbind, and participate in other reactions. We track the \emph{bound} chemical surface density separately, and label it $C_B$. We assume that at each point on the inner sphere, $C_B$ satisfies the following ODE:
\begin{align}
\frac{\partial C_B}{\partial t} = k_{on}(C^{Tot} - C_B)c_{amb} - k_{off}C_B + k_{self}C_B(C^{Tot} - C_B) + f_2,
\label{eq:coupled_cb}
\end{align}
where $k_{on}$ and $k_{off}$ are the binding and unbinding rates of $C_B$, $k_{self}$ is the rate at which $C_B$ reacts with itself, $C^{Tot}$ is the total density of binding sites at each point of the reactive zone, $f_2$ is some forcing term, and $c_{amb}$ is the concentration of the chemical $c$ at the reactive zone, \emph{i.e.}, $c_{amb} = \lf.c(\vx,t)\rt|_{\mathbb{S}_2}$. Balancing fluxes at the interface $\mathbb{S}_2$ yields the following time-varying Robin boundary condition on $c(\vx,t)$:
\begin{align}
-\nu \frac{\partial c(\vx,t)}{\partial \vn} + k_{on}(C^{Tot} - C_B)c(\vx,t) = k_{off}C_B, \vx \in \mathbb{S}_2.
\label{eq:robin_bc}
\end{align}

\subsection{A manufactured solution to the coupled problem}
The functions $f_1, f_2$, and $g$ are usually specified by platelet aggregation and coagulation models. In this article, however, we use these terms to test convergence of our numerical methods on the model problem. We again use the method of manufactured solutions, and use the terms $f_1$, and $f_2$ to make the solution hold true. The procedure is as follows. First, we set $c(\vx,t)$ to be
\begin{align}
c(\vx,t) = c(x,y,z,t) = 1 + \sin(\pi x) \cos(\pi y) \sin(\pi z) e^{-\pi t}.
\end{align}
Next, we manufacture an incompressible velocity field in the spherical shell. The velocity field $\vu(\vx,t) = (u,v,w)$ is given by
\begin{align}
\vu(\vx,t) = \phi(\vx) \sin(\pi t) [-yz,2xz,xy],
\end{align}
where $\phi(x,y,z)$ is given by
\begin{align}
\phi(\vx) = \sin\lf(\pi \|\vx\|_2^2\rt) \sin\lf(\frac{\pi}{0.3^2} \|\vx\|_2^2\rt) .
\end{align}
It is easily verified that $\vu(\vx,t)$ is incompressible in the spherical shell, and satisfies no-slip conditions on both $\mathbb{S}_1$ and $\mathbb{S}_2$; in fact, the shell domain was chosen primarily due to the ease of manufacturing this velocity field. Once $c(\vx,t)$ and $\vu(\vx,t)$ are specified, we compute the forcing term $f_1(\vx,t)$ as
\begin{align}
f_1(\vx,t) = \frac{\partial c}{\partial t} + \vu \cdot \nabla c - \nu \Delta c.
\end{align}
The boundary condition function $g(\vx,t)$ is then obtained by applying the Neumann operator $-\nu \frac{\partial}{\partial \vn}$ to $c(\vx,t)$. The function $C_B$ which is consistent with the specified $c$ is obtained by solving \eqref{eq:robin_bc} for $C_B$ to obtain
\begin{align}
C_B = \frac{-\nu \frac{\partial c}{\partial \vn} + k_{on}C^{Tot}c}{k_{on}c + k_{off}}.
\end{align}
Using this $C_B$, we compute the forcing term $f_2$ as
\begin{align}
f_2 = \frac{\partial C_B}{\partial t} - k_{on}\lf(C^{Tot}-C_B\rt)c_{amb} + k_{off}C_B - k_{self}C_B\lf(C^{Tot} - C_B\rt).
\end{align}
We compute errors against the exact $C_B$ and $c(\vx,t)$ in a spatial refinement study.

\subsection{Operator splitting for the coupled problem}

The fluid-phase chemicals $c(\vx,t)$ are coupled to the bound chemicals $C_B$ nonlinearly through the $c_{amb}$ term in \eqref{eq:coupled_cb} and the boundary conditions in \eqref{eq:robin_bc}. In order to efficiently simulate the model equations, we present a simple operator splitting algorithm based on Strang-Marchuk splitting~\cite{Strang68,Marchuk68}. For convenience, in this section we will abuse notation and use $C_B$ and $c$ to represent the numerical solutions. Given: $C_B^n$, $c^n$ (solutions at time level $n$), the algorithm is as follows:
\revtwo{
\begin{enumerate}
\item Average $C_B^n$ and $C_B^{n-1}$ to obtain $C_B^{n-1/2}$.
\item Use $C_B^{n-1/2}$ and $C_B^n$ to obtain $C_B^{n+1/2}$ by discretizing \eqref{eq:coupled_cb} with the second-order Adams-Bashforth method (AB2).
\item Compute an approximation $\tilde{C}_B^{n+1}$ to $C_B^{n+1}$ by the locally-third-order extrapolation formula $\tilde{C}_B^{n+1} = \frac{8}{3}C_B^{n+1/2} - 2C_B^n + \frac{1}{3} C_B^{n-1}$.
\item Use $\tilde{C}_B^{n+1}$ to obtain the boundary conditions for $c^{n+1}$ on $\mathbb{S}_2$ as:
\begin{align}
-\nu \frac{\partial c^{n+1}}{\partial \vn} + k_{on}\lf(C^{Tot} - \tilde{C}_B^{n+1}\rt)c^{n+1} = k_{off} \tilde{C}_B^{n+1}. 
\end{align}
\item Update $c^n$ to $c^{n+1}$ by discretizing \eqref{eq:coupled_c} with overlapped RBF-FD in space and the SBDF2 scheme in time (adding hyperviscosity implicitly in time).
\item Update $C_B^{n+1/2}$ to $C_B^{n+1}$ by computing using $c_{amb}^{n+1/2} = \lf.0.5(c^{n+1} + c^n)\rt|_{\mathbb{S}_2}$ in \eqref{eq:coupled_cb} and the AB2 discretization.
\end{enumerate}
}
It is easy to show that the splitting error is $O(\Delta t^2)$, and that the overall error in time is also $O(\Delta t^2)$. We leave the investigation of higher-order splitting schemes to future work. For the first step, we use forward Euler in place of AB2, SBDF1 in place of SBDF2, and the locally second-order extrapolation $\tilde{C}_B^{1} = 2C_B^{1/2} - C_B^0$ in place of the third-order extrapolation. We fill the initial ghost nodes using the RBF-based spatial extrapolation scheme mentioned previously.

The time-varying boundary conditions on the inner boundary $\mathbb{S}_2$ require care for efficient time-stepping. Since the boundary condition \emph{operator} itself is changing as $C_B$ changes, the entire time-stepping matrix cannot be decomposed and used in a sparse direct solve (as in the previous section). We precompute the rows of the time-stepping matrix corresponding to the domain interior $\Omega$ and the outer boundary $\mathbb{S}_1$, and \emph{append} and modify the rows corresponding to the domain boundary $\mathbb{S}_2$ every step. Since the location of $\mathbb{S}_2$ is fixed, we precompute the RBF-FD interpolation matrices and only modify the right hand sides of the linear systems corresponding to the boundary RBF-FD weights. We use the GMRES method to solve the sparse linear system for $c$~\cite{Saad86,Saad}. GMRES typically requires a good preconditioner for fast convergence. Fortunately, our problem suggests a natural preconditioner. First, we precompute an auxiliary time-stepping matrix corresponding to Neumann boundary conditions on all boundaries. We then form the sparse incomplete LU (ILU) factorization of this matrix, and use it to precondition GMRES on the problem where Robin boundary conditions are used. This allowed GMRES to converge in 2-6 iterations on average.

\revtwo{\subsection{Results}
\begin{figure}[h!]
\centering
\subfloat[]
{
	\includegraphics[scale=0.6]{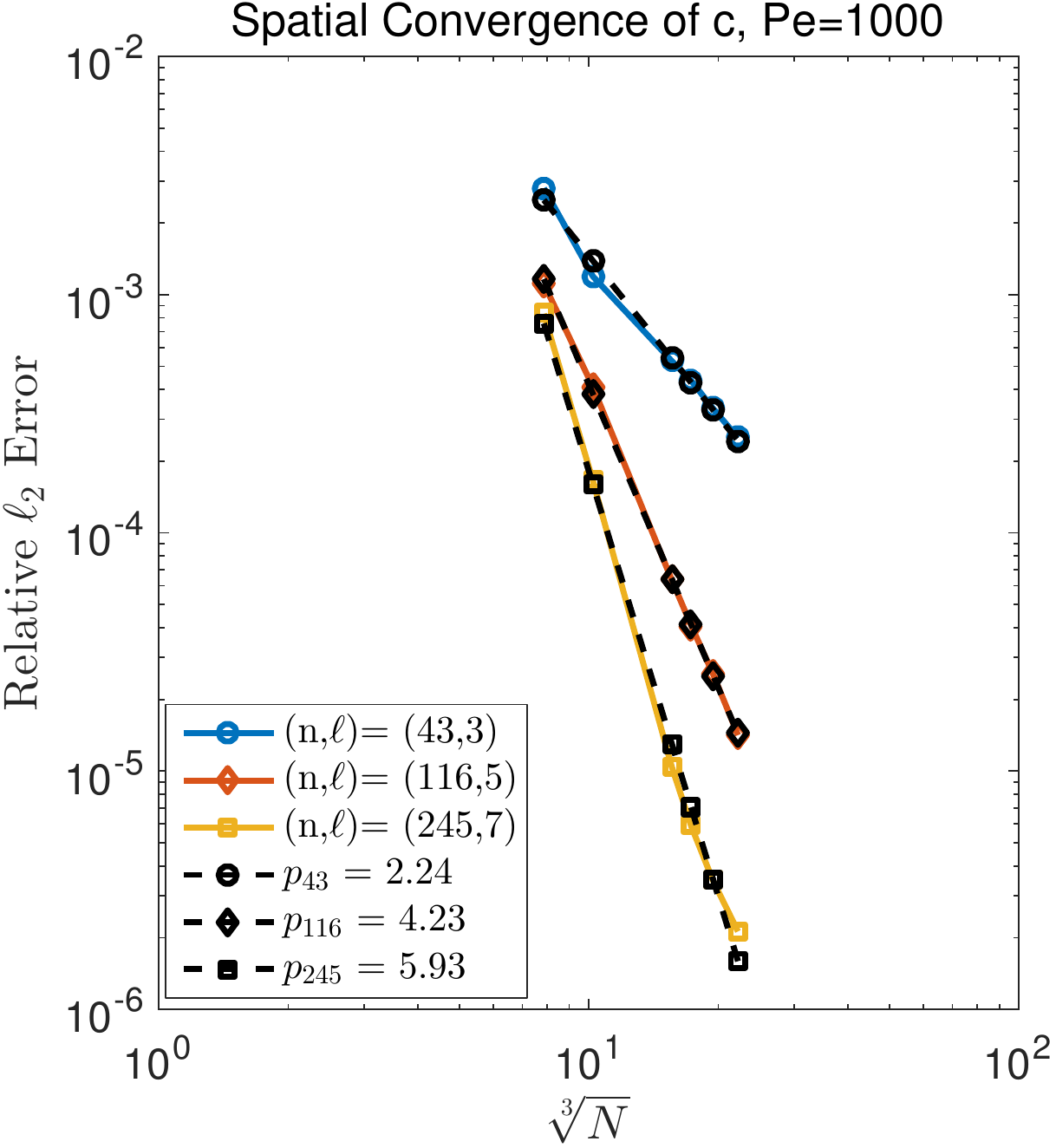}
	\label{fig:shell_11}	
}
\subfloat[]
{
	\includegraphics[scale=0.6]{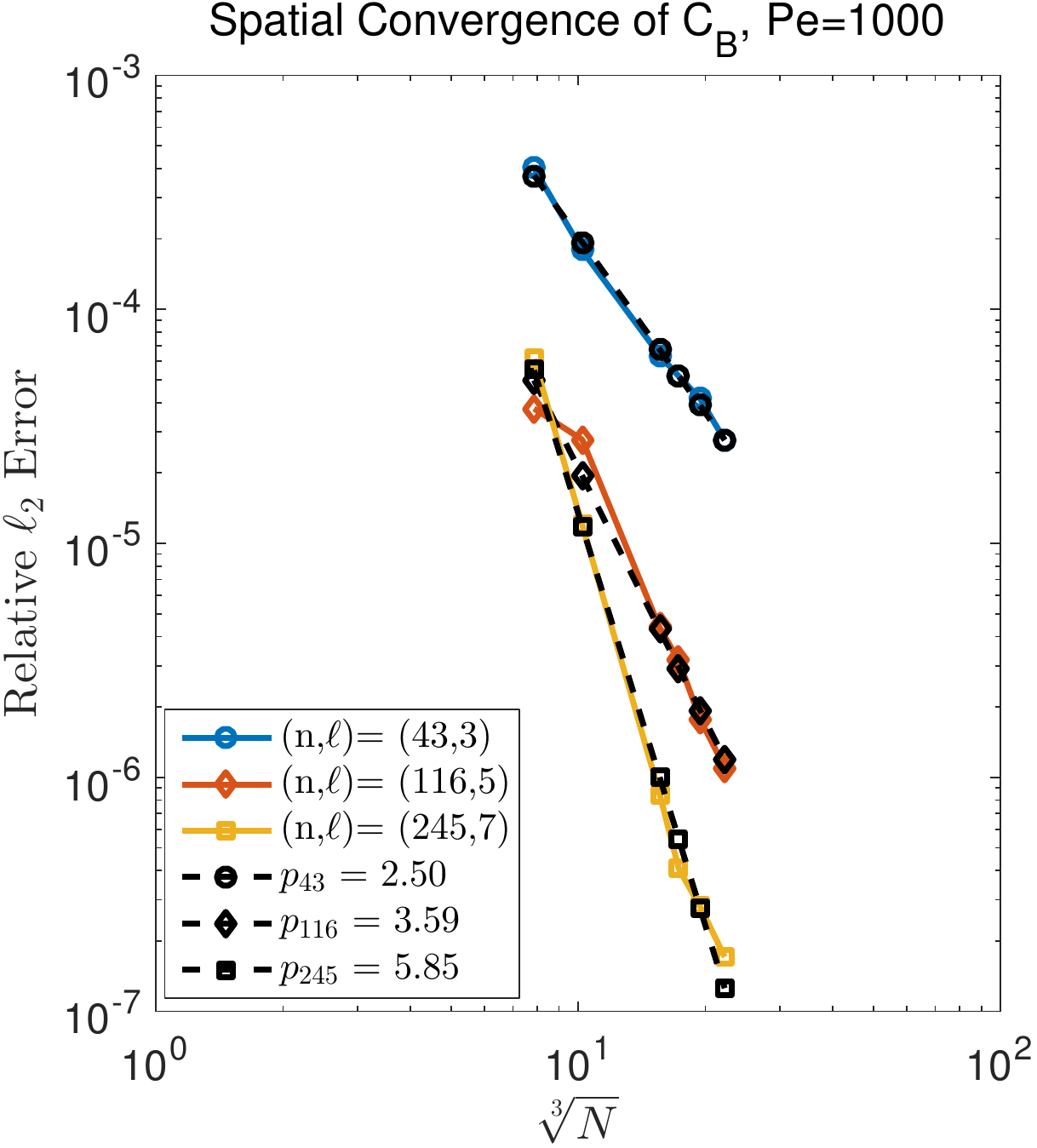}
	\label{fig:shell_21}	
}
\caption{Relative errors vs node spacing as a function of stencil size $n$ and polynomial degree $\ell$ for both $c$ (left) and $C_B$ (right).}
\label{fig:shell_conv_plots}
\end{figure}
To fully test our hyperviscosity formulation, we now set the diffusion coefficient to $\nu = 0.01$, and scale the velocity field $\vu$ to obtain a Peclet number of $\rm{Pe} = 1000$. We once again select the time-step as in the advection-diffusion test cases. We do not use any near-boundary refinement in our node sets as this may affect the CFL condition for our problem. We run convergence studies using $\xi = 2,4$ and $6$, and measure the relative errors in both the approximations to $c$ and $C_B$. The results are shown in Figure \ref{fig:shell_conv_plots}.

Figure \ref{fig:shell_11} shows that the spatial error in the numerical approximation to $c$ decreases at the rate of approximately $\xi = \ell-1$ despite the use of a second-order time-stepping scheme with time-steps close to the CFL constraint. Figure \ref{fig:shell_21} shows that the error in the approximation to $C_B$ decreases at at similar rates for $\ell=3,7$, and a slightly slower rate for $\ell=5$. At first thought, this may seem unusual since $C_B$ is governed by an ODE. However, due to our manufactured solution, $C_B$ is affected by a spatial forcing term, and thus can be expected to converge at similar rates as $c$.}

As an aside, we also remark that our discretization scheme is very well-suited to coupled problems of this sort, since our schemes automatically solve for boundary values of $c$ thereby giving us $c_{amb}$.
\section{Summary and Future Work}
\label{sec:summary}

In this article, we presented high-order numerical schemes for simulating the advection-diffusion equation using the overlapped RBF-FD method. We stabilized our methods using artificial hyperviscosity of the form $\gamma \Delta^k$. Unlike previous attempts in the RBF-FD literature, our expression for $\gamma$ was derived using a novel 1D Von-Neumann type analysis based on an explicit representation of spurious growth modes through auxiliary differential operators; the spurious growth mode was estimated numerically using a single matrix-vector multiply per differentiation matrix. We also presented an expression for $k$ based on the spectral methods literature. Using a novel ghost node formulation for IMEX time-stepping, we demonstrated high-order convergence rates on both 2D and 3D advection-diffusion problems. In addition, we demonstrated high-order convergence rates on a more complicated coupled problem in the spherical shell.

While our derived expressions for $\gamma$ are free of any tuning parameters, they require knowledge about the spectra of the differential operators being stabilized. Specifically, the quantities $\ep$ and $\eta$ must be given as an input to the algorithm. Theoretical foundations for the influence of node sets, $m$, and $\ell$ on the values of $\ep$ and $\eta$ are non-existent, forcing us to select these quantities numerically using very loose tolerances (for efficiency). We plan to explore such relationships between these parameters in future work. Also, while our article briefly explores the use of filters to preserve positivity in the numerical solution, a detailed exploration of different types of filters will be the subject of future work. Finally, our goal is to apply our hyperviscosity formulation in the setting of PDEs on time-varying domains.

\section*{Acknowledgments}
This work was supported in part on NSF grants DMS-1160432 and DMS-1521748 and NIH grant 1R01HL120728.

%\appendix
%\input{Appendix}

\section*{References}
\bibliography{article_refs_mod}

\end{document}